\baselineskip=12pt

\def\adots{\mathinner{\mskip1mu\raise1pt\hbox{.}\mskip2mu\raise4pt\hbox{.}\mskip2mu\raise7pt\vbox{\kern7pt\hbox{.}}\mskip1mu}}

\mathchardef\bfplus="062B
\mathchardef\bfminus="067B
\font\title=cmbx10 scaled\magstep5
\font\chapter=cmbx10 scaled\magstep4
\font\section=cmbx10 scaled\magstep2

\def\~#1{{\accent"7E #1}}
\def\bull{\vrule height .9ex width .8ex depth -.1ex}
\def\sqr#1#2{{\vcenter{\hrule height.#2pt \hbox{\vrule width.#2pt height#1pt \kern#1pt \vrule width.#2pt}\hrule height.#2pt}}}
\def\square{\mathchoice\sqr63\sqr63\sqr{4.2}2\sqr{1.5}2}
\def\Square{\mathop{\square}}
\def\hk#1#2{{\vcenter{\hrule height0.0pt \hbox{\vrule width0.0pt \kern#1pt \vrule width.#2pt height#1pt}\hrule height.#2pt}}}

\vskip 24pt
\centerline {\chapter Algebraically Special, Real Alpha-Geometries}
\vskip 24pt
\noindent {\section Abstract}\hfil\break
We exploit the spinor description of four-dimensional Walker geometry, and conformal rescalings of such, to describe the local geometry of four-dimensional neutral geometries with algebraically degenerate self-dual Weyl curvature and an integrable distribution of $\alpha$-planes ({\sl algebraically special real $\alpha$-geometry\/}). In particular, we determine the behaviour of Walker geometry under conformal rescaling and provide a derivation of the hyperheavenly equation from conformal rescaling formulae.
\vskip 24pt
\noindent Peter R Law. 4 Mack Place, Monroe, NY 10950, USA. prldb@member.ams.org; prldb@yahoo.com\hfil\break
Yasuo Matsushita. Section of Mathematics, School of Engineering, University of Shiga Prefecture, Hikone 522-8533, Japan. matsushita.y@e.usp.ac.jp
\vskip 24pt
\noindent 2000 MSC: 53B30, 53C27, 53C50\hfil\break
\vskip 24pt
\vfill\eject
\noindent{\section 1. Introduction}
\vskip 12pt
Throughout this paper, $(M,h)$ denotes a real, four-dimensional manifold $M$ equipped with a metric $h$ of neutral signature and will be referred to as a neutral geometry. Law (2006) contains the algebraic classification of the Weyl curvature spinors of a neutral geometry. Our aim in this paper is to shed light on those neutral geometries for which at least one of the two Weyl curvature spinors is algebraically special, i.e., admits a Weyl principal spinor (WPS), see Law (2006), of multiplicity greater than one.

An important class of examples is provided by Walker geometry, i.e., a neutral geometry $(M,g)$ that admits a parallel distribution $\cal D$ of totally null planes, see Walker (1950) and Law \& Matsushita (2008). As shown in Law \& Matsushita (2008), with its canonical choice of orientation a Walker geometry can be characterized (at least locally) by a (real) projective primed spinor field $[\pi^{A'}]$. The projective spinor field defines, at each point $p \in M$, an $\alpha$-plane
$$Z_{[\pi]} := \{\eta^A\pi^{A'}:\eta^A \in S_p\},\eqno(1.1)$$
where $S_p$ is the fibre at $p$ of the bundle (defined at least locally) of unprimed spinors and $\pi^{A'}$ is a scaled representative of $[\pi^{A'}]$. The collection of these $\alpha$-planes constitutes the distribution $\cal D$. The condition for this distribution of $\alpha$-planes (hereafter, {\sl $\alpha$-distribution\/}) to be parallel is
$$S_a := \pi_{B'}\nabla_a\pi^{B'} = 0,\eqno(1.2)$$
where $\pi^{A'}$ is any local scaled representative (LSR) of $[\pi^{A'}]$, i.e., any local primed spinor field whose projective class coincides with $[\pi^{A'}]$ (at each point). The condition for the $\alpha$-distribution to be merely integrable is
$$\pi_{B'}\pi^{A'}\nabla_{AA'}\pi^{B'} = 0.\eqno(1.3)$$
Clearly (1.2) entails (1.3). Equation (1.2) is studied in Law \& Matsushita (2008), where it is shown that a solution $[\pi^{A'}]$ of (1.2) is a WPS of multiplicity at least two, and preliminary results for equation (1.3) presented in Law (2008). We will denote a Walker geometry by $(M,g,Z_{[\pi]})$. For simplicity, we will phrase our discussion as if any $(M,h)$ admits projective spinor bundles, which presumes that $(M,h)$ is ${\bf SO}^\bfplus$-orientable (where ${\bf SO^\bfplus}$ refers to the identity component of {\bf O(2,2)}); otherwise the results are valid locally.

The Generalized Goldberg-Sachs Theorem (GGST) for neutral geometry (see Law 2008 (6.2.17)) indicates an intimate connexion between solutions of (1.3) and multiple WPSs (real or complex; note that while the spinor spaces for neutral signature are real, solutions of spinorial equations may be complex valued, complex WPSs are just a particular instance). We denote by $(M,h,[\pi^{A'}])$ a neutral geometry for which $[\pi^{A'}]$ is a solution of (1.3) and refer to such as an $\alpha$-{\sl geometry}. Every solution of (1.3) is automatically a WPS (Law 2008 (6.2.9)). We call $(M,h,[\pi^{A'}])$ an {\sl algebraically special\/} (AS) $\alpha$-geometry when $[\pi^{A'}]$ is a multiple WPS (of the SD Weyl curvature spinor $\tilde\Psi_{A'B'C'D'}$). We leave untouched the possibility of neutral geometries $(M,h)$ with a WPS $[\pi^{A'}]$, of multiplicity $p > 1$, which is not a solution of (1.3). By GGST, for such a neutral geometry:
$$\underbrace{\pi^{A'}\ldots\pi^{C'}}_{5-p}\nabla^{DD'}\tilde\Psi_{A'B'C'D'} \not= 0.\eqno(1.4)$$
Note that if a complex $[\pi^{A'}]$ is is a WPS, then its conjugate is too and its multiplicity is at most two.

We refer to neutral geometries admitting real (complex) solutions of (1.3) as, respectively, real (complex) $\alpha$-geometries. The cases of real and complex WPSs will be treated separately as the geometric interpretation of (1.3) differs in the two situations; in this paper, we treat real $\alpha$-geometries and real multiple WPSs, and we may omit the qualifier `real' when it is not needed for emphasis.

Consider a Walker geometry $(M,g,Z_{[\pi]})$ and its integrable $\alpha$-distribution $Z = Z_{[\pi]}$, defined as in (1.1). An $\alpha$-distribution is clearly a conformally invariant notion and integrability is a differential-topological property. Hence, in the neutral geometry $(M,h)$, where $h := \Omega^2g$ for some smooth ${\bf R}^+$-valued function $\Omega$, the distribution $Z$ retains its character as the integrable $\alpha$-distribution $Z_{[\pi]}$. Of course, $[\pi^{A'}]$ is also a multiple WPS for $(M,h)$ since the Weyl curvature spinors are invariant under conformal rescalings of the metric. Thus, conformal rescalings of Walker geometries generate further examples of real AS$\alpha$-geometries, which raises the question: is every real AS$\alpha$-geometry $(M,h,[\pi^{A'}])$ conformally Walker, i.e., is $h$ a conformal rescaling of a metric $g$, with $(M,g,Z_{[\pi]})$ Walker? In \S 2, we review the geometry of conformal rescalings in the four-dimensional, neutral-signature context. In \S 3, we deduce the behaviour of Walker geometry under conformal rescalings (3.1--26) and, further, show that every real AS$\alpha$-geometry is {\sl locally} conformal to a Walker geometry (3.27), thereby obtaining a local characterization of real AS$\alpha$-geometry by applying conformal rescaling formulae (2.8) to the curvature spinors for Walker geometry derived in Law \& Matsushita (2008).

In the context of complex general relativity, Pleba\'nski and co-workers, see Pleba\'nski \& Robinson (1976, 1977), Finley \& Pleba\'nski (1976), and Boyer et al. (1980), studied the condition of algebraically special SD Weyl curvature and showed that the vacuum Einstein equations reduce to a single nonlinear PDE, which they called the {\sl hyperheavenly equation}. In \S 3, we show (3.30--3.70) how the Ricci curvature of conformally rescaled Walker geometry leads to (the real neutral-signature version of) the hyperheavenly equation (3.71--73). Subsequent to our work, we learnt that Pleba\'nski \& R\'ozga (1984) had presented a version of our (3.27) (their Lemma IV, though without reference to Walker geometry), but apparently this result was not exploited and no connexion with Walker's (1950) work noted. Indeed, in retrospect, it seems remarkable that Walker (1950) appears never to have been exploited in complex general relativity.

Recently, however, Chudecki \& Przanowski (2008b) presented Walker geometry as a special case of the hyperheavenly formalism applied to neutral geometry. Our perspective, on the other hand, is that Walker geometry, characterized by a parallel totally null distribution, is a more fundamental geometric setting. As shown in Law (2008), the Einstein condition for Walker geometries $(M,g,Z_{[\pi]})$ corresponds to the so-called nonexpanding hyperheavenly equation. In this paper, we in effect show that the so-called expanding hyperheavenly equation corresponds to the Einstein condition for neutral metrics which are locally conformally Walker. But our approach, with its emphasis on conformal rescaling, the spinor formalism of Penrose \& Rindler (1984), and Walker geometry, differs from, and we feel is more natural than, the direct application of the hyperheavenly formalism to the study of neutral geometry. Rather than just presenting a new derivation of the hyperheavenly equation, we regard our work as a study of the behaviour of Walker geometry under conformal rescalings, viz., both (3.1) and (3.27), which illuminates the local geometry of algebraically special real $\alpha$-geometries.

In \S 4, we touch on some related null geometry.

As in Law \& Matsushita (2008) and Law (2008), we employ the abstract index formalism and notation of Penrose \& Rindler (1984) for indices (italic indices are `abstract' indices while bold upright indices are `concrete'). See Law (2006), Law \& Matsushita (2008) Appendix Two, and Law (2008) for brief accounts of the adaptation of the two-component spinor formalism of Penrose \& Rindler (1984) to the context of neutral signature. $S$ will denote the bundle of unprimed spinors over $M$ (or over an open set $U$ if $S$ only exists locally) while $S'$ denotes the bundle of primed spinors.
\vskip 24pt
\noindent{\section 2. Conformal Rescalings of the Metric}
\vskip 12pt
Penrose \& Rindler (1984), \S 5.6, describe how, in the context of Lorentzian four-dimensional geometry, a spinor structure is related to a conformal class of metrics. This correspondence is easily adapted to the case of neutral signature. Let $(M,h)$ be a neutral geometry with, at least locally, a spinor structure. Let $\hat h_{ab} := \Omega^2h_{ab}$, where $\Omega:M \to {\bf R}^+$ is smooth. The corresponding spinor structure for $(M,\hat h)$ consists of the same spinor bundles but with skew scalar products given by
$$\hat\epsilon_{AB} := \theta\epsilon_{AB} \hskip 1in \hat\epsilon_{A'B'} := \tilde\theta\epsilon_{A'B'},\eqno(2.1)$$
where $\theta$ and $\tilde\theta$ are smooth functions $M \to {\bf R}^+$ and
$$\theta\tilde\theta = \Omega^2.\eqno(2.2)$$
In the usual manner, the Levi-Civita connection induces unique connections (denoted simply by $\hat\nabla_{AA'}$) on the unprimed and primed spinor bundles with respect to which $\hat\epsilon_{AB}$ and $\hat\epsilon_{A'B'}$, respectively, are parallel. By the neutral-signature analogue of Penrose \& Rindler (1984) (4.4.23),
$$\hat\nabla_{AA'}\xi^C = \nabla_{AA'}\xi^C + \Theta_{AA'B}{}^C\xi^B \hskip 1in \hat\nabla_{AA'}\xi^{C'} = \nabla_{AA'}\xi^{C'} + \tilde\Theta_{AA'B'}{}^{C'}\xi^{B'},\eqno(2.3)$$
where $\Theta_{aB}{}^C$ and $\tilde\Theta_{aB'}{}^{C'}$ are, in principle, independent quantities. Since, however, both $\hat\nabla_{AA'}$ and $\nabla_{AA'}$ are torsion free, one deduces, following Penrose \& Rindler (1984), pp. 216--217, that $\Theta_{aB}{}^C$ and $\tilde\Theta_{aB'}{}^{C'}$ take the form
$$\Theta_{aB}{}^C = \Lambda_a\epsilon_B{}^C + \Upsilon_{A'B}\epsilon_A{}^C \hskip 1in \tilde\Theta_{aB'}{}^{C'} = \tilde\Lambda_a\epsilon_{B'}{}^{C'} + \tilde\Upsilon_{AB'}\epsilon_{A'}{}^{C'},$$
but with
$$\Lambda_a + \tilde\Lambda_a = 0 \hskip 1.25in \Upsilon_a = \tilde\Upsilon_a,$$
in place of Penrose \& Rindler (1984) (4.4.44 \& 46). Hence
$$\Theta_{aB}{}^C = \Lambda_a\epsilon_B{}^C + \Upsilon_{A'B}\epsilon_A{}^C \hskip 1in \tilde\Theta_{aB'}{}^{C'} = -\Lambda_a\epsilon_{B'}{}^{C'} + \Upsilon_{AB'}\epsilon_{A'}{}^{C'}.\eqno(2.4)$$
Now,
$$\eqalign{0 &= \hat\nabla_a\hat\epsilon_{BC}\cr
&= \hat\nabla_a(\theta\epsilon_{BC})\cr
&= \nabla_a(\theta\epsilon_{BC}) - \Theta_{aB}{}^D\epsilon_{DC}\theta - \Theta_{aC}{}^D\epsilon_{BD}\theta\cr
&= \epsilon_{BC}\nabla_a(\theta) - \epsilon_{BC}\Theta_{aD}{}^D\theta\cr}$$
whence
$$\nabla_a(\ln\theta) = \Theta_{aD}{}^D = 2\Lambda_a + \Upsilon_a\qquad\hbox{and similarly}\qquad \nabla_a(\ln\tilde\theta) = \tilde\Theta_{aD'}{}^{D'} = -2\Lambda_a + \Upsilon_a.\eqno(2.5)$$
With $\mu:M \to {\bf R}^+$ smooth, write $\theta = \mu\Omega$, whence $\tilde\theta = \Omega/\mu$. Then (2.5) is equivalent to
$$\Upsilon_a = \nabla_a(\ln\Omega) \hskip 1.25in \Lambda_a = \nabla_a\left(\ln\sqrt\mu\right).\eqno(2.6)$$
In the case of Lorentzian signature, the fact that the primed spin space is the complex conjugate of the unprimed spin space forces $\tilde\theta$ to be the complex conjugate of $\theta$, whence $\mu\overline\mu = 1$. Provided $\mu$ is required to be real valued, then $\mu = \pm 1$. The choice of negative sign is rejected as being discontinuous with the identity scaling, and $\mu \equiv 1$ results in the formulae of Penrose \& Rindler (1984) \S 5.6. The independence of the primed and unprimed spin spaces for neutral signature is less restrictive and does allow nontrivial choices of $\mu$. Clearly, $\Lambda_a = 0$ is equivalent to choosing $\mu$ constant. Whether nontrivial choices of $\mu$ are of interest in neutral geometry will not be pursued here; hereafter, we suppose $\mu \equiv 1$, whence unprimed and primed spinors are treated alike under conformal rescaling and $\hat\nabla_a$ is determined in terms of $\nabla_a$ and $\Upsilon_a = \nabla_a(\ln\Omega)$ exactly as in Penrose \& Rindler (1984) \S 5.6. In particular,
$$\eqalignno{\hat\nabla_{AA'}\chi^{P\ldots S'\ldots}_{B\ldots F'\ldots} &= \nabla_{AA'}\chi^{P\ldots S'\ldots}_{B\ldots F'\ldots} - \Upsilon_{BA'}\chi^{P\ldots S'\ldots}_{A\ldots F'\ldots} - \cdots - \Upsilon_{AF'}\chi^{P\ldots S'\ldots}_{B\ldots A'\ldots} - \cdots\cr
&+ \epsilon_A{}^P\Upsilon_{XA'}\chi^{X\ldots S'\ldots}_{B\ldots F'\ldots} + \cdots + \epsilon_{A'}{}^{S'}\Upsilon_{AX'}\chi^{P\ldots X'\ldots}_{B\ldots F'\ldots} + \cdots&(2.7)\cr}$$
Our conventions for curvature are detailed in Law \& Matsushita (2008), Appendix One; they agree with those of Penrose \& Rindler (1984), except that our Ricci tensor (whence Ricci scalar curvature) is the negative of theirs. But we take the Ricci spinor and the scalar $\Lambda$ to be unaffected by this different convention (one simply inserts an additional negative sign in Penrose \& Rindler (1984)(4.6.20--23)). The curvature spinors of the metric $\hat h_{ab}$ are then related to the curvature spinors of the metric $h_{ab}$ as in Penrose \& Rindler (1984):
$$\displaylines{\hat{\tilde\Psi}_{A'B'C'D'} = \tilde\Psi_{A'B'C'D'} \hskip 1.25in \hat\Psi_{ABCD} = \Psi_{ABCD}\cr
\hfill\hat\Phi_{ABA'B'} = \Phi_{ABA'B'} + \Upsilon_{A(A'}\Upsilon_{B')B} - \nabla_{A(A'}\Upsilon_{B')B}\hfill\llap(2.8)\cr
\hat\Lambda = \Omega^{-2}\left[\Lambda + {1 \over 4}\nabla^a\Upsilon_a + {1 \over 4}\Upsilon^a\Upsilon_a\right] = \Omega^{-2}\left[\Lambda + {1 \over 4}\Omega^{-1}\Square\Omega\right].\cr}$$
It will prove useful to define
$$\omega := \ln\Omega.\eqno(2.9)$$
Notation for spin coefficients for neutral geometry (effectively a suitable definition of the priming operation to replace Penrose \& Rindler (1984) (4.5.17)) was introduced in Law (2008). The notation of Law (2008) combined with (2.9) yields, for the components of $\Upsilon_a$ with respect to spin bases $\{o^A,\iota^A\}$ and $\{o^{A'},\iota^{A'}\}$,
$$\Upsilon_{00'} = D\omega \qquad \Upsilon_{01'} = \delta\omega \qquad \Upsilon_{10'} = \triangle\omega \qquad \Upsilon_{11'} = D'\omega,\eqno(2.10)$$
which equations in effect define the operators $D$, $\delta$, $\triangle$ and $D'$.
If the spin basis $\{o^A,\iota^A\}$ is rescaled according to
$$\hat o^A := \Omega^{v_0}o^A \hskip 1in \hat\iota^A := \Omega^{v_1}\iota^A,\eqno(2.11)$$
then
$$\hat o_A = \Omega^{v_0 + 1}o_A \hskip 1in \hat\iota_A = \Omega^{v_1 + 1}\iota_A,\eqno(2.12)$$
with
$$\hat\chi := \hat\iota^A\hat o_A = \Omega^{v_0+v_1+1}\chi,\eqno(2.13)$$
where $\chi :=  \iota^Ao_A$. The spin basis $\{o^{A'},\iota^{A'}\}$ may be independently rescaled as
$$\hat o^{A'} := \Omega^{w_0}o^{A'} \hskip 1in \hat\iota^{A'} := \Omega^{w_1}\iota^{A'},\eqno(2.14)$$
whence
$$\hat o_{A'} = \Omega^{w_0 + 1}o_{A'} \hskip 1in \hat\iota_{A'} = \Omega^{v_1 + 1}\iota_{A'},\eqno(2.15)$$
with
$$\hat{\tilde\chi} := \hat\iota^{A'}\hat o_{A'} = \Omega^{w_0+w_1+1}\tilde\chi,\eqno(2.16)$$
where $\tilde\chi :=  \iota^{A'}o_{A'}$. With these choices, as $\hat\nabla_a = \nabla_a$ on functions, then on functions
$$\hat D = \Omega^{v_0+w_0}D \qquad \hat \delta = \Omega^{v_0+w_1}\delta \qquad \hat\triangle = \Omega^{v_1+w_0}\triangle \qquad \hat D' = \Omega^{v_1+w_1}D'.\eqno(2.17)$$
 
The spin coefficients of $\hat\nabla_a$ with respect to the rescaled bases may be derived from Law (2008) (2.7--8), which are the neutral analogues of Penrose \& Rindler (1984) (4.5.21). One obtains the following neutral analogue of Penrose \& Rindler (1984) (5.6.25):
\vskip 12pt
$$\vcenter{\offinterlineskip
\halign{&\vrule#&\strut\ #\ \cr
\noalign{\hrule}
height3pt&\omit&&\omit&&\omit&&\omit&\cr
&\hfil$\hat\epsilon$\hfil&&\hfil$\hat\kappa$\hfil&&\hfil$\hat\tau'$\hfil&&\hfil$\hat\gamma'$\hfil&\cr
height3pt&\omit&&\omit&&\omit&&\omit&\cr
\noalign{\hrule}
height3pt&\omit&&\omit&&\omit&&\omit&\cr
&\hfil$\hat\alpha$\hfil&&\hfil$\hat\rho$\hfil&&\hfil$\hat\sigma'$\hfil&&\hfil$\hat\beta'$\hfil&\cr
height3pt&\omit&&\omit&&\omit&&\omit&\cr
\noalign{\hrule}
height3pt&\omit&&\omit&&\omit&&\omit&\cr
&\hfil$\hat\beta$\hfil&&\hfil$\hat\sigma$\hfil&&\hfil$\hat\rho'$\hfil&&\hfil$\hat\alpha'$\hfil&\cr
height3pt&\omit&&\omit&&\omit&&\omit&\cr
\noalign{\hrule}
height3pt&\omit&&\omit&&\omit&&\omit&\cr
&\hfil$\hat\gamma$\hfil&&\hfil$\hat\tau$\hfil&&\hfil$\hat\kappa'$\hfil&&\hfil$\hat\epsilon'$\hfil&\cr
height3pt&\omit&&\omit&&\omit&&\omit&\cr
\noalign{\hrule}}}
=
\vcenter{\offinterlineskip
\halign{\strut\ #\ \cr
\noalign{\vskip 3pt}
\hfil$\Omega^{w_0+v_1}.$\hfil\cr
\noalign{\vskip 3pt}
\hfil$\Omega^{w_0+v_1}.$\hfil\cr
\noalign{\vskip 3pt}
\hfil$\Omega^{v_0+w_1}.$\hfil\cr
\noalign{\vskip 3pt}
\hfil$\Omega^{v_0+w_1}.$\hfil\cr
\noalign{\vskip 3pt}}}
\vcenter{\offinterlineskip
\halign{&\vrule#&\strut\ #\ \cr
\noalign{\hrule}
height3pt&\omit&&\omit&&\omit&&\omit&\cr
&\hfil$[\epsilon+(v_0+1)D\omega]\Sigma$\hfil&&\hfil$\kappa\Sigma^2$\hfil&&\hfil$\tau'-\triangle\omega$\hfil&&\hfil$[\gamma'+v_1D\omega]\Sigma$\hfil&\cr
height3pt&\omit&&\omit&&\omit&&\omit&\cr
\noalign{\hrule}
height3pt&\omit&&\omit&&\omit&&\omit&\cr
&\hfil$\alpha+v_0\triangle\omega$\hfil&&\hfil$[\rho+D\omega]\Sigma$\hfil&&\hfil$\sigma'\Sigma^{-1}$\hfil&&\hfil$\beta'-(v_1+1)\triangle\omega$\hfil&\cr
height3pt&\omit&&\omit&&\omit&&\omit&\cr
\noalign{\hrule}
height3pt&\omit&&\omit&&\omit&&\omit&\cr
&\hfil$\beta+(v_0+1)\delta\omega$\hfil&&\hfil$\sigma\Sigma$\hfil&&\hfil$[\rho'+D'\omega]\Sigma^{-1}$\hfil&&\hfil$\alpha'-v_1\delta\omega$\hfil&\cr
height3pt&\omit&&\omit&&\omit&&\omit&\cr
\noalign{\hrule}
height3pt&\omit&&\omit&&\omit&&\omit&\cr
&\hfil$[\gamma+v_0D'\omega]\Sigma^{-1}$\hfil&&\hfil$\tau+\delta\omega$\hfil&&\hfil$\kappa'\Sigma^{-2}$\hfil&&\hfil$[\epsilon'+(v_1+1)D'\omega]\Sigma^{-1}$\hfil&\cr
height3pt&\omit&&\omit&&\omit&&\omit&\cr
\noalign{\hrule}}}\eqno(2.18)$$
\vskip 24pt
$$\vcenter{\offinterlineskip
\halign{&\vrule#&\strut\ #\ \cr
\noalign{\hrule}
height3pt&\omit&&\omit&&\omit&&\omit&\cr
&\hfil$\hat{\tilde\epsilon}$\hfil&&\hfil$\hat{\tilde\kappa}$\hfil&&\hfil$\skew{-4}\hat{\tilde\tau'}$\hfil&&\hfil$\skew{-4}\hat{\tilde\gamma'}$\hfil&\cr
height3pt&\omit&&\omit&&\omit&&\omit&\cr
\noalign{\hrule}
height3pt&\omit&&\omit&&\omit&&\omit&\cr
&\hfil$\hat{\tilde\alpha}$\hfil&&\hfil$\hat{\tilde\rho}$\hfil&&\hfil$\skew{-4}\hat{\tilde\sigma'}$\hfil&&\hfil$\skew{-4}\hat{\tilde\beta'}$\hfil&\cr
height3pt&\omit&&\omit&&\omit&&\omit&\cr
\noalign{\hrule}
height3pt&\omit&&\omit&&\omit&&\omit&\cr
&\hfil$\hat{\tilde\beta}$\hfil&&\hfil$\hat{\tilde\sigma}$\hfil&&\hfil$\skew{-4}\hat{\tilde\rho'}$\hfil&&\hfil$\skew{-4}\hat{\tilde\alpha'}$\hfil&\cr
height3pt&\omit&&\omit&&\omit&&\omit&\cr
\noalign{\hrule}
height3pt&\omit&&\omit&&\omit&&\omit&\cr
&\hfil$\hat{\tilde\gamma}$\hfil&&\hfil$\hat{\tilde\tau}$\hfil&&\hfil$\skew{-4}\hat{\tilde\kappa'}$\hfil&&\hfil$\skew{-4}\hat{\tilde\epsilon'}$\hfil&\cr
height3pt&\omit&&\omit&&\omit&&\omit&\cr
\noalign{\hrule}}}
=
\vcenter{\offinterlineskip
\halign{\strut\ #\ \cr
\noalign{\vskip 3pt}
\hfil$\Omega^{v_0+w_1}.$\hfil\cr
\noalign{\vskip 3pt}
\hfil$\Omega^{v_0+w_1}.$\hfil\cr
\noalign{\vskip 3pt}
\hfil$\Omega^{w_0+v_1}.$\hfil\cr
\noalign{\vskip 3pt}
\hfil$\Omega^{w_0+v_1}.$\hfil\cr
\noalign{\vskip 3pt}}}
\vcenter{\offinterlineskip
\halign{&\vrule#&\strut\ #\ \cr
\noalign{\hrule}
height3pt&\omit&&\omit&&\omit&&\omit&\cr
&\hfil$[\tilde\epsilon+(w_0+1)D\omega]\tilde\Sigma$\hfil&&\hfil$\tilde\kappa\tilde\Sigma^2$\hfil&&\hfil$\tilde\tau'-\delta\omega$\hfil&&\hfil$[\tilde\gamma'+w_1D\omega]\tilde\Sigma$\hfil&\cr
height3pt&\omit&&\omit&&\omit&&\omit&\cr
\noalign{\hrule}
height3pt&\omit&&\omit&&\omit&&\omit&\cr
&\hfil$\tilde\alpha+w_0\delta\omega$\hfil&&\hfil$[\tilde\rho+D\omega]\tilde\Sigma$\hfil&&\hfil$\tilde\sigma'\tilde\Sigma^{-1}$\hfil&&\hfil$\tilde\beta'-(w_1+1)\delta\omega$\hfil&\cr
height3pt&\omit&&\omit&&\omit&&\omit&\cr
\noalign{\hrule}
height3pt&\omit&&\omit&&\omit&&\omit&\cr
&\hfil$\tilde\beta+(w_0+1)\triangle\omega$\hfil&&\hfil$\tilde\sigma\tilde\Sigma$\hfil&&\hfil$[\tilde\rho'+D'\omega]\tilde\Sigma^{-1}$\hfil&&\hfil$\tilde\alpha'-w_1\triangle\omega$\hfil&\cr
height3pt&\omit&&\omit&&\omit&&\omit&\cr
\noalign{\hrule}
height3pt&\omit&&\omit&&\omit&&\omit&\cr
&\hfil$[\tilde\gamma+w_0D'\omega]\tilde\Sigma^{-1}$\hfil&&\hfil$\tilde\tau+\triangle\omega$\hfil&&\hfil$\tilde\kappa'\tilde\Sigma^{-2}$\hfil&&\hfil$[\tilde\epsilon'+(w_1+1)D'\omega]\tilde\Sigma^{-1}$\hfil&\cr
height3pt&\omit&&\omit&&\omit&&\omit&\cr
\noalign{\hrule}}}\eqno(2.19)$$
\vskip 12pt
where
$$\Sigma := \Omega^{v_0-v_1} \hskip 1.25in \tilde\Sigma := \Omega^{w_0-w_1},\eqno(2.20)$$
and the factors of $\Omega$ to the right of the equality signs in (2.18--19) multiply each entry in the corresponding row.
\vskip 24pt
\noindent {\section 3. Algebraically Special Real Alpha-Geometries}
\vskip 12pt
Let $(M,g,Z_{[\pi]})$ be a Walker geometry. As discussed in \S 2, the $\alpha$-distribution $Z_{[\pi]}$ of the Walker geometry retains its character as an integrable (though not necessarily parallel) $\alpha$-distribution with respect to the conformal class of metrics $[g]$ and may unambiguously be denoted by $Z_{[\pi]}$ with respect to $[g]$. 

Consider the neutral geometry $(M,\hat g)$, $\hat g := \Omega^2g$ with $\Omega:M \to {\bf R}^+$ smooth. First note that as long as LSRs of $[\pi^{A'}]$ are treated as conformal densities/invariants, then their status as LSRs is valid for the conformal class of metrics to which $g$ belongs. With $\hat\pi^{A'} := \Omega^p\pi^{A'}$, ($p=0$ allowed), a simple calculation yields:
$$\hat\pi_{A'}\hat\pi^{B'}\hat\nabla_b\hat\pi^{A'} = \Omega^{3p+1}\pi_{A'}\pi^{B'}\hat\nabla_b\pi^{A'} = \Omega^{3p+1}\pi_{A'}\pi^{B'}[\nabla_b\pi^{A'} + \epsilon_{B'}{}^{A'}\Upsilon_{BX'}\pi^{X'}] = \Omega^{3p+1}\pi_{A'}\pi^{B'}\nabla_b\pi^{A'},$$
i.e., as expected $\hat\pi^{A'}$ satisfies (1.3) for $(M,\Omega^2g)$ iff $\pi^{A'}$ solves (1.3) for $(M,g)$.

Since $[\pi^{A'}]$ is WPS of multiplicity at least two (Law \& Matsushita 2008, 2.5) for $(M,g)$, by (2.8) it is so for $(M,\Omega^2g)$, whence, as asserted in the Introduction, $(M,\Omega^2g)$ is indeed a real AS$\alpha$-geometry (and it is easy to check using (2.8) that $\hat\pi^{A'}\hat\pi^{B'}\hat\pi^{C'}\hat\nabla^{DD'}\hat{\tilde\Psi}_{A'B'C'D'} = 0$). It is natural to ask when the $\alpha$-distribution $Z_{[\pi]}$ is parallel in $(M,\hat g)$, i.e., when is $(M,\Omega^2g,[\pi^{A'}])$ itself Walker? A simple computation like the last yields the following result.
\vskip 24pt
\noindent {\bf 3.1 Lemma}\hfil\break
If $(M,g,Z_{[\pi]})$ is Walker, then $(M,\Omega^2g,[\pi^{A'}])$ is Walker iff $\Upsilon_{BX'}\pi^{X'} = 0$ for some, whence any, LSR of $[\pi^{A'}]$, i.e., iff $\Upsilon^a \in Z_{[\pi]}$, in other words, iff $\Omega$ is constant on $\alpha$-surfaces. When $(M,\Omega^2g,Z_{[\pi]})$ is Walker, for any LSR, $\hat\nabla_b\pi^{A'} = \nabla_b\pi^{A'}$, whence $(M,g,Z_{[\pi]})$ admits a parallel LSR iff $(M,\Omega^2g,Z_{[\pi]})$ does (see Law \& Matsushita 2008, \S 3, for parallel LSRs in Walker geometry).

Proof. If $\hat\pi^{A'} := \Omega^p\pi^{A'}$, then
$$\hat\pi_{A'}\hat\nabla_b\hat\pi^{A'} = \Omega^{2p+1}\pi_{A'}[\nabla_b\pi^{A'} + \epsilon_{B'}{}^{A'}\Upsilon_{BX'}\pi^{X'}] = \Omega^{2p+1}\pi_{B'}\Upsilon_{BX'}\pi^{X'},$$
whence $\hat\pi^{A'}$ solves (1.2) iff $\Upsilon_{BX'}\pi^{X'} = 0$. Thus, $(M,\Omega^2g,[\pi^{A'}])$ is Walker iff $\Upsilon^a \in Z_{[\pi]}$. When this condition holds, then clearly $\hat\nabla_b\pi^{A'} = \nabla_b\pi^{A'}$.\bull
\vskip 24pt
Suppose $(u,v,x,y)$ are Walker coordinates for $(M,g,Z_{[\pi]})$ on some domain, i.e., they are Frobenius coordinates for the $\alpha$-distribution and the metric components with respect to these coordinates take the form 
$$\left(g_{\bf ab}\right) = \pmatrix{{\bf 0}_2&{\bf 1}_2\cr {\bf 1}_2&W\cr} \hskip 1in W = \pmatrix{a&c\cr c&b\cr} =: \left(W^{\bf AB}\right),\eqno(3.2)$$
where $a$, $b$, and $c$ are functions of $(u,v,x,y)$, see Walker (1950), Law \& Matsushita (2008) 2.3. Whether $(M,\hat g,[\pi^{A'}])$ is Walker or not, the coordinates $(u,v,x,y)$ are Frobenius coordinates for the $\alpha$-distribution irrespective of the metric and still provide useful coordinates for $(M,\hat g,[\pi^{A'}])$, the metric $\hat g_{ab}$ having components
$$\left(\hat g_{\bf ab}\right) = \Omega^2\pmatrix{{\bf 0}_2&{\bf 1}_2\cr {\bf 1}_2&W\cr}.\eqno(3.3)$$
We will call such coordinates {\sl conformal Walker coordinates} for $(M,\hat g,[\pi^{A'}])$.

There are at least two simple and natural choices of null tetrads related to conformal Walker coordinates. First, since the coordinates $x$ and $y$ label distinct $\alpha$-surfaces in the chart domain, for any LSR $\pi^{A'}$, one can write
$$dx = \mu_A\pi_{A'} \hskip 1in dy = \nu_A\pi_{A'},\qquad \nu^A\mu_A \not= 0.\eqno(3.4)$$
As in Law \& Matsushita (2008) 2.4, given (conformal) Walker coordinates $(u,v,x,y)$, one can choose a unique (up to sign)  LSR $\pi^{A'}$ so that (after possibly interchanging $u$ and $v$ and $x$ and $y$) Law \& Matsushita (2008) (2.8) holds:
$$dx = \alpha_A\pi_{A'} \qquad dy = \beta_A\pi_{A'} \qquad \partial_u = \alpha^A\pi^{A'} \qquad \partial_v = \beta^A\pi^{A'},\eqno(3.5)$$
where $\{\alpha^A,\beta^A\}$ is a spin frame ($\beta^A\alpha_A = 1$). We will refer to such (conformal) Walker coordinates as {\sl oriented} because one can choose an atlas of such iff the $\alpha$-distribution is orientable, $M$ itself being naturally oriented by the atlas of all Walker coordinates (see Law \& Matsushita 2008 \S 1). The Walker null tetrad for $(M,g,Z_{[\pi]})$ introduced in Law \& Matsushita (2008) (2.11) is
$$\vcenter{\openup1\jot \halign{$\hfil#$&&${}#\hfil$&\qquad$\hfil#$\cr
\ell_a &= dx = \alpha_A\pi_{A'} & \tilde m_a &= dy = \beta_A\pi_{A'}\cr
n_a &= \displaystyle du + {a \over 2}dx + {c \over 2}dy = \beta_A\xi_{A'} & m_a &= \displaystyle-(dv + {c \over 2}dx + {b \over 2}dy) = \alpha_A\xi_{A'}\cr
\ell^a &= \partial_u = \alpha^A\pi^{A'} & \tilde m^a &= \partial_v = \beta^A\pi^{A'}\cr
n^a &= \displaystyle-{a \over 2}\partial_u - {c \over 2}\partial_v + \partial_x = \beta^A\xi^{A'} & m^a &= \displaystyle{c \over 2}\partial_u + {b \over 2}\partial_v - \partial_y = \alpha^A\xi^{A'},\cr}}\eqno(3.6)$$
\vskip 6pt
where $\{\pi^{A'},\xi^{A'}\}$ is indeed a spin frame for the Walker geometry (these two spin frames corresponding to the Walker null tetrad are called Walker spin frames). If one takes $L_a := \ell_a = dx$ and $\tilde M_a = \tilde m_a = dy$, observe that $dx \wedge dy = \epsilon_{AB}\pi_{A'}\pi_{B'} = \Omega^{-1}\hat\epsilon_{AB}\pi_{A'}\pi_{B'}$, which suggests taking $\hat \pi_{A'} := \Omega^{-1/2}\pi_{A'}$ as an LSR so that 
$$dx \wedge dy = 2\ell_{[a}\tilde m_{b]} = \epsilon_{AB}\pi_{A'}\pi_{B'} = \hat\epsilon_{AB}\hat\pi_{A'}\hat\pi_{B'} = 2L_{[a}\tilde M_{b]}.\eqno(3.7)$$
To maintain $L_a = \ell_a$ and $\tilde M_a = \tilde m_a$, write
$$\ell_a = \alpha_A\pi_{A'} = \hat\alpha_A\hat\pi_{A'} = L_a \hskip 1in \tilde m_A = \beta_A\pi_{A'} = \hat\beta_A\hat\pi_{A'} = \tilde M_a.\eqno(3.8)$$
If one further defines $N_a := \Omega^2n_a$, $M_a := \Omega^2m_a$, and $\hat\xi_{A'} := \Omega^{3/2}\xi_{A'}$, then
$$\vcenter{\openup1\jot \halign{$\hfil#$&&${}#\hfil$&\qquad$\hfil#$\cr
\hat\pi_{A'} &= \Omega^{-1/2}\pi_{A'} & \hat\xi_{A'} &= \Omega^{3/2}\xi_{A'} & \hat\alpha_A &= \Omega^{1/2}\alpha_A & \hat\beta_A &= \Omega^{1/2}\beta_A\cr
\hat\pi^{A'} &= \Omega^{-3/2}\pi^{A'} & \hat\xi^{A'} &= \Omega^{1/2}\xi^{A'} & \hat\alpha^A &= \Omega^{-1/2}\alpha^A & \hat\beta^A &= \Omega^{-1/2}\beta^A,\cr}}\eqno(3.9)$$
where, of course, hatted objects have indices raised and lowered with respect to the geometry of $(M,\hat g)$. Furthermore,
$$\vcenter{\openup1\jot \halign{$\hfil#$&&${}#\hfil$&\qquad$\hfil#$\cr
L_a &= \ell_a & \tilde M_a &= \tilde m_a & N_a &= \Omega^2n_a & M_a &= \Omega^2m_a\cr
L^a &= \Omega^{-2}\ell^a & \tilde M^a &= \Omega^{-2}\tilde m^a & N^a &= n^a & M^a &= m^a,\cr}}\eqno(3.10)$$
and $\{L^a,N^a,M^a,\tilde M^a\}$ is a null tetrad for the geometry $(M,\hat g)$ ($\hat g_{ab} = 2L_{(a}N_{b)} - 2M_{(a}\tilde M_{b)}$). Moreover, $\{\hat\alpha^a,\hat\beta^A\}$ and $\{\hat\pi^{A'},\hat\xi^{A'}\}$ are spin frames in the geometry of $(M,\hat g)$ and, indeed, the spin frames associated to the null tetrad $\{L^a,N^a,M^a,\tilde M^a\}$, so that
$$L^a = \hat\alpha^A\hat\pi^{A'} \qquad N^a = \hat\beta^A\hat\xi^{A'} \qquad M^a = \hat\alpha^A\hat\xi^{A'} \qquad \tilde M^a = \hat\beta^A\hat\pi^{A'}.\eqno(3.11)$$
This choice of spin frames for $(M,\hat g)$ is of the form (2.11--15), with 
$$v_0 = v_1 = -{1 \over 2} \qquad w_0 = -{3 \over 2} \qquad w_1 = {1 \over 2},\eqno(3.12)$$
whence the spin coefficients for $(M,\hat g)$ with respect to these spin frames can be obtained directly from (2.18--19) using the formulae for the spin coefficients for the Walker geometry $(M,g)$ with respect to the Walker spin frames given in Law (2008) (5.6). The spin frames $\{\hat\alpha^a,\hat\beta^A\}$ and $\{\hat\pi^{A'},\hat\xi^{A'}\}$ are also convenient for taking components of the curvature spinors in (2.8) as the right hand side can be evaluated in terms of the components of the curvature spinors for $(M,g)$ with respect to the Walker spin frames, see Law \& Matsushita (2008) \S 2, and components of $\Upsilon_a$. In particular,
\vskip 12pt
$$\vcenter{\offinterlineskip
\halign{&\vrule#&\strut\ #\ \cr
\noalign{\hrule}
height3pt&\omit&&\omit&&\omit&&\omit&\cr
&\hfil$\hat\epsilon$\hfil&&\hfil$\hat\kappa$\hfil&&\hfil$\hat\tau'$\hfil&&\hfil$\hat\gamma'$\hfil&\cr
height3pt&\omit&&\omit&&\omit&&\omit&\cr
\noalign{\hrule}
height3pt&\omit&&\omit&&\omit&&\omit&\cr
&\hfil$\hat\alpha$\hfil&&\hfil$\hat\rho$\hfil&&\hfil$\hat\sigma'$\hfil&&\hfil$\hat\beta'$\hfil&\cr
height3pt&\omit&&\omit&&\omit&&\omit&\cr
\noalign{\hrule}
height3pt&\omit&&\omit&&\omit&&\omit&\cr
&\hfil$\hat\beta$\hfil&&\hfil$\hat\sigma$\hfil&&\hfil$\hat\rho'$\hfil&&\hfil$\hat\alpha'$\hfil&\cr
height3pt&\omit&&\omit&&\omit&&\omit&\cr
\noalign{\hrule}
height3pt&\omit&&\omit&&\omit&&\omit&\cr
&\hfil$\hat\gamma$\hfil&&\hfil$\hat\tau$\hfil&&\hfil$\hat\kappa'$\hfil&&\hfil$\hat\epsilon'$\hfil&\cr
height3pt&\omit&&\omit&&\omit&&\omit&\cr
\noalign{\hrule}}}
=
\vcenter{\offinterlineskip
\halign{\strut\ #\ \cr
\noalign{\vskip 3pt}
\hfil$\Omega^{-2}.$\hfil\cr
\noalign{\vskip 3pt}
\hfil$\Omega^{-2}.$\hfil\cr
\noalign{\vskip 3pt}
\hfil$\Omega^{0}.$\hfil\cr
\noalign{\vskip 3pt}
\hfil$\Omega^{0}.$\hfil\cr
\noalign{\vskip 3pt}}}
\vcenter{\offinterlineskip
\halign{&\vrule#&\strut\ #\ \cr
\noalign{\hrule}
height3pt&\omit&&\omit&&\omit&&\omit&\cr
&\hfil${1 \over 2}D\omega$\hfil&&\hfil$0$\hfil&&\hfil$-\triangle\omega$\hfil&&\hfil$-{1\over 2}D\omega$\hfil&\cr
height3pt&\omit&&\omit&&\omit&&\omit&\cr
\noalign{\hrule}
height3pt&\omit&&\omit&&\omit&&\omit&\cr
&\hfil$-{1 \over 2}\triangle\omega$\hfil&&\hfil$D\omega$\hfil&&\hfil$0$\hfil&&\hfil$-{1 \over 2}\triangle\omega$\hfil&\cr
height3pt&\omit&&\omit&&\omit&&\omit&\cr
\noalign{\hrule}
height3pt&\omit&&\omit&&\omit&&\omit&\cr
&\hfil$\beta+{1 \over 2}\delta\omega$\hfil&&\hfil$\sigma$\hfil&&\hfil$\rho'+D'\omega$\hfil&&\hfil$\alpha'+{1 \over 2}\delta\omega$\hfil&\cr
height3pt&\omit&&\omit&&\omit&&\omit&\cr
\noalign{\hrule}
height3pt&\omit&&\omit&&\omit&&\omit&\cr
&\hfil$\gamma-{1\over 2}D'\omega$\hfil&&\hfil$\tau+\delta\omega$\hfil&&\hfil$\kappa'$\hfil&&\hfil$\epsilon'+{1\over 2}D'\omega$\hfil&\cr
height3pt&\omit&&\omit&&\omit&&\omit&\cr
\noalign{\hrule}}}\eqno(3.13)$$
\vskip 24pt
$$\vcenter{\offinterlineskip
\halign{&\vrule#&\strut\ #\ \cr
\noalign{\hrule}
height3pt&\omit&&\omit&&\omit&&\omit&\cr
&\hfil$\hat{\tilde\epsilon}$\hfil&&\hfil$\hat{\tilde\kappa}$\hfil&&\hfil$\skew{-4}\hat{\tilde\tau'}$\hfil&&\hfil$\skew{-4}\hat{\tilde\gamma'}$\hfil&\cr
height3pt&\omit&&\omit&&\omit&&\omit&\cr
\noalign{\hrule}
height3pt&\omit&&\omit&&\omit&&\omit&\cr
&\hfil$\hat{\tilde\alpha}$\hfil&&\hfil$\hat{\tilde\rho}$\hfil&&\hfil$\skew{-4}\hat{\tilde\sigma'}$\hfil&&\hfil$\skew{-4}\hat{\tilde\beta'}$\hfil&\cr
height3pt&\omit&&\omit&&\omit&&\omit&\cr
\noalign{\hrule}
height3pt&\omit&&\omit&&\omit&&\omit&\cr
&\hfil$\hat{\tilde\beta}$\hfil&&\hfil$\hat{\tilde\sigma}$\hfil&&\hfil$\skew{-4}\hat{\tilde\rho'}$\hfil&&\hfil$\skew{-4}\hat{\tilde\alpha'}$\hfil&\cr
height3pt&\omit&&\omit&&\omit&&\omit&\cr
\noalign{\hrule}
height3pt&\omit&&\omit&&\omit&&\omit&\cr
&\hfil$\hat{\tilde\gamma}$\hfil&&\hfil$\hat{\tilde\tau}$\hfil&&\hfil$\skew{-4}\hat{\tilde\kappa'}$\hfil&&\hfil$\skew{-4}\hat{\tilde\epsilon'}$\hfil&\cr
height3pt&\omit&&\omit&&\omit&&\omit&\cr
\noalign{\hrule}}}
=
\vcenter{\offinterlineskip
\halign{\strut\ #\ \cr
\noalign{\vskip 3pt}
\hfil$\Omega^{0}.$\hfil\cr
\noalign{\vskip 3pt}
\hfil$\Omega^{0}.$\hfil\cr
\noalign{\vskip 3pt}
\hfil$\Omega^{-2}.$\hfil\cr
\noalign{\vskip 3pt}
\hfil$\Omega^{-2}.$\hfil\cr
\noalign{\vskip 3pt}}}
\vcenter{\offinterlineskip
\halign{&\vrule#&\strut\ #\ \cr
\noalign{\hrule}
height3pt&\omit&&\omit&&\omit&&\omit&\cr
&\hfil$-{1\over 2}\Omega^{-2}D\omega$\hfil&&\hfil$0$\hfil&&\hfil$-\delta\omega$\hfil&&\hfil${1\over 2}\Omega^{-2}D\omega$\hfil&\cr
height3pt&\omit&&\omit&&\omit&&\omit&\cr
\noalign{\hrule}
height3pt&\omit&&\omit&&\omit&&\omit&\cr
&\hfil$\tilde\alpha-{3 \over 2}\delta\omega$\hfil&&\hfil$\Omega^{-2}D\omega$\hfil&&\hfil$\tilde\sigma'\Omega^2$\hfil&&\hfil$\tilde\beta'-{3 \over 2}\delta\omega$\hfil&\cr
height3pt&\omit&&\omit&&\omit&&\omit&\cr
\noalign{\hrule}
height3pt&\omit&&\omit&&\omit&&\omit&\cr
&\hfil$-{1\over 2}\triangle\omega$\hfil&&\hfil$0$\hfil&&\hfil$\Omega^2D'\omega$\hfil&&\hfil$\tilde\alpha'-{1 \over 2}\triangle\omega$\hfil&\cr
height3pt&\omit&&\omit&&\omit&&\omit&\cr
\noalign{\hrule}
height3pt&\omit&&\omit&&\omit&&\omit&\cr
&\hfil$[\tilde\gamma-{3 \over 2}D'\omega]\Omega^2$\hfil&&\hfil$\triangle\omega$\hfil&&\hfil$\tilde\kappa'\Omega^4$\hfil&&\hfil$[\tilde\epsilon'+{3 \over 2}D'\omega]\Omega^2$\hfil&\cr
height3pt&\omit&&\omit&&\omit&&\omit&\cr
\noalign{\hrule}}}\eqno(3.14)$$
\vskip 12pt
\noindent where formulae for the nonzero spin coefficients of $(M,g)$ with respect to the Walker spin frames in terms of $a$, $b$, and $c$ can be substituted from Law (2008) (5.6).

Alternatively, one may instead require $L^a = \ell^a$ and $\tilde M^a = \tilde m^a$ so as to maintain this local frame for $Z_{[\pi]}$. By similar reasoning as above, since now
$$2L^{[a}\tilde M^{b]} = 2\ell^{[a}\tilde m^{b]} = \epsilon^{AB}\pi^{A'}\pi^{B'} = \Omega\hat\epsilon^{AB}\pi^{A'}\pi^{B'},$$
one is led to define
$$\hat\pi^{A'} := \Omega^{1/2}\pi^{A'} \qquad \hat\xi^{A'} := \Omega^{-3/2}\xi^{A'} \qquad \hat\alpha^A := \Omega^{-1/2}\alpha^A \qquad \hat\beta^A := \Omega^{-1/2}\beta^A,\eqno(3.15)$$
whence
$$\hat\pi_{A'} = \Omega^{3/2}\pi_{A'} \qquad \hat\xi_{A'} = \Omega^{-1/2}\xi_{A'} \qquad \hat\alpha_A = \Omega^{1/2}\alpha_A \qquad \hat\beta_A = \Omega^{1/2}\beta_A.\eqno(3.16)$$
Now
$$L^a := \hat\alpha^A\hat\pi^{A'} = \ell^a \qquad N^a:= \hat\beta^A\hat\xi^{A'} = \Omega^{-2}n^a \qquad M^a := \hat\alpha^A\hat\xi^{A'} = \Omega^{-2}m^a \qquad \tilde M^a := \hat\beta^A\hat\pi^{A'} = \tilde m^a,\eqno(3.17)$$
and
$$L_a = \Omega^2\ell_a \qquad N_a = n_a \qquad M_a = m_a \qquad \tilde M_a = \Omega^2\tilde m_a,\eqno(3.18)$$
so
$$\hat g_{ab} = 2L_{(a}N_{b)} - 2M_{(a}\tilde M_{b)} = \Omega^2(2\ell_{(a}n_{b)} - 2m_{(a}\tilde m_{b)}) = \Omega^2g_{ab} .\eqno(3.19)$$
Of course, $\{\hat\alpha^A,\hat\beta^A\}$ and $\{\hat\pi^{A'},\hat\xi^{A'}\}$ are the spin frames associated to the null tetrad $\{L^a,N^a,M^a,\tilde M^a\}$ for $(M,\hat g)$. This alternative choice of spin frames is again of the form (2.11--15), but with
$$v_0 = v_1 = -{1 \over 2} \qquad w_0 = {1 \over 2} \qquad w_1 = -{3 \over 2},\eqno(3.20)$$
whence
\vskip 12pt
$$\vcenter{\offinterlineskip
\halign{&\vrule#&\strut\ #\ \cr
\noalign{\hrule}
height3pt&\omit&&\omit&&\omit&&\omit&\cr
&\hfil$\hat\epsilon$\hfil&&\hfil$\hat\kappa$\hfil&&\hfil$\hat\tau'$\hfil&&\hfil$\hat\gamma'$\hfil&\cr
height3pt&\omit&&\omit&&\omit&&\omit&\cr
\noalign{\hrule}
height3pt&\omit&&\omit&&\omit&&\omit&\cr
&\hfil$\hat\alpha$\hfil&&\hfil$\hat\rho$\hfil&&\hfil$\hat\sigma'$\hfil&&\hfil$\hat\beta'$\hfil&\cr
height3pt&\omit&&\omit&&\omit&&\omit&\cr
\noalign{\hrule}
height3pt&\omit&&\omit&&\omit&&\omit&\cr
&\hfil$\hat\beta$\hfil&&\hfil$\hat\sigma$\hfil&&\hfil$\hat\rho'$\hfil&&\hfil$\hat\alpha'$\hfil&\cr
height3pt&\omit&&\omit&&\omit&&\omit&\cr
\noalign{\hrule}
height3pt&\omit&&\omit&&\omit&&\omit&\cr
&\hfil$\hat\gamma$\hfil&&\hfil$\hat\tau$\hfil&&\hfil$\hat\kappa'$\hfil&&\hfil$\hat\epsilon'$\hfil&\cr
height3pt&\omit&&\omit&&\omit&&\omit&\cr
\noalign{\hrule}}}
=
\vcenter{\offinterlineskip
\halign{\strut\ #\ \cr
\noalign{\vskip 3pt}
\hfil$\Omega^{0}.$\hfil\cr
\noalign{\vskip 3pt}
\hfil$\Omega^{0}.$\hfil\cr
\noalign{\vskip 3pt}
\hfil$\Omega^{-2}.$\hfil\cr
\noalign{\vskip 3pt}
\hfil$\Omega^{-2}.$\hfil\cr
\noalign{\vskip 3pt}}}
\vcenter{\offinterlineskip
\halign{&\vrule#&\strut\ #\ \cr
\noalign{\hrule}
height3pt&\omit&&\omit&&\omit&&\omit&\cr
&\hfil${1 \over 2}D\omega$\hfil&&\hfil$0$\hfil&&\hfil$-\triangle\omega$\hfil&&\hfil$-{1\over 2}D\omega$\hfil&\cr
height3pt&\omit&&\omit&&\omit&&\omit&\cr
\noalign{\hrule}
height3pt&\omit&&\omit&&\omit&&\omit&\cr
&\hfil$-{1 \over 2}\triangle\omega$\hfil&&\hfil$D\omega$\hfil&&\hfil$0$\hfil&&\hfil$-{1 \over 2}\triangle\omega$\hfil&\cr
height3pt&\omit&&\omit&&\omit&&\omit&\cr
\noalign{\hrule}
height3pt&\omit&&\omit&&\omit&&\omit&\cr
&\hfil$\beta+{1 \over 2}\delta\omega$\hfil&&\hfil$\sigma$\hfil&&\hfil$\rho'+D'\omega$\hfil&&\hfil$\alpha'+{1 \over 2}\delta\omega$\hfil&\cr
height3pt&\omit&&\omit&&\omit&&\omit&\cr
\noalign{\hrule}
height3pt&\omit&&\omit&&\omit&&\omit&\cr
&\hfil$\gamma-{1\over 2}D'\omega$\hfil&&\hfil$\tau+\delta\omega$\hfil&&\hfil$\kappa'$\hfil&&\hfil$\epsilon'+{1\over 2}D'\omega$\hfil&\cr
height3pt&\omit&&\omit&&\omit&&\omit&\cr
\noalign{\hrule}}}\eqno(3.21)$$
\vskip 24pt
$$\vcenter{\offinterlineskip
\halign{&\vrule#&\strut\ #\ \cr
\noalign{\hrule}
height3pt&\omit&&\omit&&\omit&&\omit&\cr
&\hfil$\hat{\tilde\epsilon}$\hfil&&\hfil$\hat{\tilde\kappa}$\hfil&&\hfil$\skew{-4}\hat{\tilde\tau'}$\hfil&&\hfil$\skew{-4}\hat{\tilde\gamma'}$\hfil&\cr
height3pt&\omit&&\omit&&\omit&&\omit&\cr
\noalign{\hrule}
height3pt&\omit&&\omit&&\omit&&\omit&\cr
&\hfil$\hat{\tilde\alpha}$\hfil&&\hfil$\hat{\tilde\rho}$\hfil&&\hfil$\skew{-4}\hat{\tilde\sigma'}$\hfil&&\hfil$\skew{-4}\hat{\tilde\beta'}$\hfil&\cr
height3pt&\omit&&\omit&&\omit&&\omit&\cr
\noalign{\hrule}
height3pt&\omit&&\omit&&\omit&&\omit&\cr
&\hfil$\hat{\tilde\beta}$\hfil&&\hfil$\hat{\tilde\sigma}$\hfil&&\hfil$\skew{-4}\hat{\tilde\rho'}$\hfil&&\hfil$\skew{-4}\hat{\tilde\alpha'}$\hfil&\cr
height3pt&\omit&&\omit&&\omit&&\omit&\cr
\noalign{\hrule}
height3pt&\omit&&\omit&&\omit&&\omit&\cr
&\hfil$\hat{\tilde\gamma}$\hfil&&\hfil$\hat{\tilde\tau}$\hfil&&\hfil$\skew{-4}\hat{\tilde\kappa'}$\hfil&&\hfil$\skew{-4}\hat{\tilde\epsilon'}$\hfil&\cr
height3pt&\omit&&\omit&&\omit&&\omit&\cr
\noalign{\hrule}}}
=
\vcenter{\offinterlineskip
\halign{\strut\ #\ \cr
\noalign{\vskip 3pt}
\hfil$\Omega^{-2}.$\hfil\cr
\noalign{\vskip 3pt}
\hfil$\Omega^{-2}.$\hfil\cr
\noalign{\vskip 3pt}
\hfil$\Omega^{0}.$\hfil\cr
\noalign{\vskip 3pt}
\hfil$\Omega^{0}.$\hfil\cr
\noalign{\vskip 3pt}}}
\vcenter{\offinterlineskip
\halign{&\vrule#&\strut\ #\ \cr
\noalign{\hrule}
height3pt&\omit&&\omit&&\omit&&\omit&\cr
&\hfil${3\over 2}\Omega^{2}D\omega$\hfil&&\hfil$0$\hfil&&\hfil$-\delta\omega$\hfil&&\hfil$-{3\over 2}\Omega^{2}D\omega$\hfil&\cr
height3pt&\omit&&\omit&&\omit&&\omit&\cr
\noalign{\hrule}
height3pt&\omit&&\omit&&\omit&&\omit&\cr
&\hfil$\tilde\alpha+{1 \over 2}\delta\omega$\hfil&&\hfil$\Omega^{2}D\omega$\hfil&&\hfil$\tilde\sigma'\Omega^{-2}$\hfil&&\hfil$\tilde\beta'+{1 \over 2}\delta\omega$\hfil&\cr
height3pt&\omit&&\omit&&\omit&&\omit&\cr
\noalign{\hrule}
height3pt&\omit&&\omit&&\omit&&\omit&\cr
&\hfil${3\over 2}\triangle\omega$\hfil&&\hfil$0$\hfil&&\hfil$\Omega^{-2}D'\omega$\hfil&&\hfil${3 \over 2}\triangle\omega$\hfil&\cr
height3pt&\omit&&\omit&&\omit&&\omit&\cr
\noalign{\hrule}
height3pt&\omit&&\omit&&\omit&&\omit&\cr
&\hfil$[\tilde\gamma+{1 \over 2}D'\omega]\Omega^{-2}$\hfil&&\hfil$\triangle\omega$\hfil&&\hfil$\tilde\kappa'\Omega^{-4}$\hfil&&\hfil$[\tilde\epsilon'-{1 \over 2}D'\omega]\Omega^{-2}$\hfil&\cr
height3pt&\omit&&\omit&&\omit&&\omit&\cr
\noalign{\hrule}}}\eqno(3.22)$$
\vskip 12pt
\noindent where again the nonzero spin coefficients on the right-hand sides are those for $(M,g)$ with respect to the Walker spin frames and are given in Law (2008) (5.6) in terms of $a$, $b$, and $c$.

In particular, both (3.14) and (3.22) confirm that $\hat{\tilde\kappa} = \hat{\tilde\sigma} = 0$, which are the conditions for the $\alpha$-distribution to be integrable (see Law 2008 (6.2.4)) and that $(M,\hat g,[\pi^{A'}])$ is Walker iff $\hat{\tilde\rho} = \hat{\tilde\tau} = 0$ (see Law 2008 (6.2.40)), i.e., iff $D\omega = \triangle\omega = 0$ as required by (3.1). 

When $(M,\hat g,Z_{[\pi]})$ is itself Walker, (3.13-14) and (3.21--22) simplify considerably. One also readily observes, for example, that, employing either set of spin frames (3.9) or (3.15--16), by (2.8)
$$\hat\Phi_{ABA'B'}\hat\pi^{A'}\hat\pi^{B'} = \Omega^k[\Phi_{ABA'B'}\pi^{A'}\pi^{B'} + \Upsilon_{AA'}\pi^{A'}\Upsilon_{BB'}\pi^{B'} - \pi^{A'}\pi^{B'}\nabla_{AA'}\Upsilon_{BB'}],$$
for an appropriate value of $k$. The first term on the right-hand side vanishes by Law \& Matsushita (2008) 2.5; the second term directly from (3.1); and, upon writing $\Upsilon_{BB'} = \kappa_B\pi_{B'}$ by virtue of (3.1), the third term is seen to vanish by (1.3). Thus, $\hat\Phi_{ABA'B'}\hat\pi^{A'}\hat\pi^{B'}$ vanishes, confirming that $[\pi^{A'}]$ is a Ricci principal spinor (RPS), as it must be in a Walker geometry (Law \& Matsushita 2008, 2.5).

Clearly, even when $(M,\hat g,Z_{[\pi]})$ is Walker, the metric components of $\hat g_{ab} = \Omega^2g_{ab}$ with respect to $(u,v,x,y)$ are not of the form (3.2), i.e, $(u,v,x,y)$ are not Walker coordinates for $(M,\hat g,Z_{[\pi]})$, and the null tetrads and spin frames of (3.9--10) and (3.15--18) are not {\sl Walker} null tetrads nor {\sl Walker} spin frames. But it is easy to construct Walker coordinates for $(M,\hat g,Z_{[\pi]})$ from $(u,v,x,y)$ if desired. Since the coordinates $x$ and $y$ label distinct $\alpha$-surfaces in the chart domain, for any LSR $\pi^{A'}$, as usual one can write
$$dx = \mu_A\pi_{A'} \hskip 1in dy = \nu_A\pi_{A'},\qquad \nu^A\mu_A \not= 0.\eqno(3.23)$$
Now define 
$$\hat U^a := \hat g^{ab}\mu_B\pi_{B'} = \Omega^{-2}g^{ab}\mu_B\pi_{B'} =: \Omega^{-2}U^a \hskip .5in \hat V^a := \hat g^{ab}\nu_B\pi_{B'} = \Omega^{-2}g^{ab}\nu_B\pi_{B'} =: \Omega^{-2}V^a.\eqno(3.24)$$
Following Law \& Matsushita (2008) 2.3, one can choose solutions $\hat u$ and $\hat v$ of the system of equations
$$\hat U\hat u = 1 = \hat V\hat v \hskip 1in \hat U\hat v = 0 = \hat V \hat u,$$
i.e.,
$$U\hat u = \Omega^2 = V\hat v \hskip 1in U\hat v = 0 = V \hat u.$$
The coordinates $u$ and $v$ are constructed in Law \& Matsushita (2008) 2.3 so that $U = \partial_u$ and $V = \partial_v$. By (3.1), $\Omega$ is constant on $\alpha$-surfaces, i.e., a function of $x$ and $y$ only, when expressed in terms of the Walker coordinates $(u,v,x,y)$. Hence, one can choose
$$\hat u :=\Omega^2u \qquad \hat v := \Omega^2v \qquad \hat x := x \qquad \hat y := y,\eqno(3.25)$$
and $(\hat u,\hat v,\hat x, \hat y)$ are Walker coordinates for $(M,\hat g,[\pi^{A'}])$. Indeed,
\vskip 6pt
$$(\hat\partial_{\bf j}) := \pmatrix{\partial_{\hat u}&\partial_{\hat v}&\partial_{\hat x}& \partial_{\hat y}\cr} = \pmatrix{\partial_{u}&\partial_{v}&\partial_{x}&\partial_{y}\cr}.\pmatrix{\Omega^{-2}&0&-2\Omega^{-3}\Omega_x\hat u&-2\Omega^{-3}\Omega_y\hat u\cr 0&\Omega^{-2}&-2\Omega^{-3}\Omega_x\hat v&-2\Omega^{-3}\Omega_y\hat v\cr 0&0&1&0\cr 0&0&0&1\cr} =: (\partial_{\bf i})(J^{\bf i}{}_{\bf j})$$
\vskip 6pt
whence the components of $\hat g_{ab}$ with respect to $(\hat u,\hat v,\hat x,\hat y)$ are
$${^\tau\! J}.\Omega^2\pmatrix{{\bf 0}_2&{\bf 1}_2\cr {\bf 1}_2&W\cr}.J = \pmatrix{{\bf 0}_2&{\bf 1}_2\cr {\bf 1}_2&\hat W\cr}$$
where
$$\hat W := \Omega^2\pmatrix{a - 4\Omega^{-3}\Omega_x\hat u&c - 2\Omega^{-3}(\Omega_x\hat v + \Omega_y\hat u)\cr c-2\Omega^{-3}(\Omega_x\hat v + \Omega_y\hat u)&b-4\Omega^{-3}\Omega_y\hat v \cr} =: \pmatrix{\hat a&\hat c\cr \hat c&\hat b\cr}.$$
If $(u,v,x,y)$ are oriented Walker coordinates for $(M,g,[\pi^{A'}])$ and the LSR $\pi^{A'}$ is chosen as in (3.5), then
$$dx \wedge dy = \epsilon_{AB}\pi_{A'}\pi_{B'} = \hat\epsilon_{AB}\hat\pi_{A'}\hat\pi_{B'}$$
if one puts $\hat\pi_{A'} := \Omega^{-1/2}\pi_{A'}$; equivalently,
$$\partial_{\hat u} \wedge \partial_{\hat v} = \Omega^{-4}\partial_u \wedge \partial_v = \Omega^{-4}\epsilon^{AB}\pi^{A'}\pi^{B'} = \hat\epsilon^{AB}\hat\pi^{A'}\hat\pi^{B'},$$
where $\hat\pi^{A'} = \Omega^{-3/2}\pi^{A'}$. Thus, $(\hat u,\hat v,\hat x,\hat y)$ are oriented Walker coordinates for $(M,\hat g,[\pi^{A'}])$ and $\hat\pi^{A'}$ is an LSR for these coordinates as in (3.5). If $\{L^a,N^a,M^a,\tilde M^a\}$ denotes the Walker null tetrad for $(M,\hat g,[\pi^{A'}])$ and the coordinates $(\hat u,\hat v,\hat x,\hat y)$, then $L_a = d\hat x = dx = \ell_a = \alpha_A\pi_{A'} = \hat\alpha_A\hat\pi_{A'}$ if $\hat\alpha_A := \Omega^{1/2}\alpha_A$. Similarly, $\tilde M_a = d\hat y = dy = \tilde m_a$ suggests putting $\hat\beta_A := \Omega^{1/2}\beta_A$. From (3.6),
$$\displaylines{N_a = d\hat u + {\hat a \over 2}d\hat x + {\hat c\over 2}d\hat y = \Omega^2n_a + \Omega(\Omega_yu-\Omega_xv)\tilde m_a = \hat\beta_A\bigl[\Omega^{3/2}\xi_{A'} + \Omega^{1/2}(\Omega_yu-\Omega_xv)\pi_{A'}\bigr];\cr
M_a = -\left(d\hat v + {\hat c\over2}d\hat x + {\hat b \over 2}d\hat y\right) = \Omega^2m_a + \Omega(\Omega_yu-\Omega_xv)\ell_a = \hat\alpha_A\bigl[\Omega^{3/2}\xi_{A'} + \Omega^{1/2}(\Omega_yu-\Omega_xv)\pi_{A'}\bigr].\cr}$$
One can now identify the Walker spin frames $\{\hat\alpha^A,\hat\beta^A\}$ and $\{\hat\pi^{A'},\hat\xi^{A'}\}$ for $(M,\hat g,[\pi^{A'}])$ and the coordinates $(\hat u,\hat v,\hat x,\hat y)$:
$$\vcenter{\openup2\jot \halign{$\hfil#$&&${}#\hfil$&\qquad$\hfil#$\cr
\hat\alpha^A &:= \Omega^{-1/2}\alpha^A & \hat\beta^A &:= \Omega^{-1/2}\beta^A\cr
\hat\pi^{A'} &:= \Omega^{-3/2}\pi^{A'} & \hat\xi^{A'} &:= \Omega^{1/2}\xi^{A'} + \Omega^{-1/2}(\Omega_yu-\Omega_xv)\pi^{A'}.\cr}}\eqno(3.26)$$
The primed Walker spin frames for $(M,\hat g,[\pi^{A'}])$ and $(M,g,[\pi^{A'}])$ are not related as in (2.14) so (2.19) is inapplicable. Of course, the spin coefficients for the Walker spin frames for $(M,\hat g,[\pi^{A'}])$ may be obtained directly from Law (2008) (5.6).

After these preliminaries on Walker geometry, we turn to the consideration of arbitrary real AS$\alpha$-geometries.
\vskip 24pt
\noindent {\bf 3.27 Proposition}\hfil\break
A real $\alpha$-geometry $(M,h,[\pi^{A'}])$ is algebraically special iff locally conformally Walker, i.e., iff for each $p \in M$, there exists a neighbourhood $U_p$ and a metric $g$ such that $(U_p,g,[\pi^{A'}])$ is Walker and $h = \Omega^2g$ on $U_p$ for some $\Omega:U_p \to {\bf R}^+$.

Proof. We seek a neighbourhood $U_p$ and function $\chi:U_p \to {\bf R}^+$ such that $(U_p,\chi^2h,Z_{[\pi]})$ is Walker, i.e., such that for any, whence every, LSR $\pi^{A'}$, (1.2) holds with respect to the metric $\hat h := \chi^2h$. Hence, one requires a $\chi$ such that
$$0 = \pi_{A'}\hat\nabla_b\pi^{A'} = \pi_{A'}[\nabla_b\pi^{A'} + \epsilon_{B'}{}^{A'}\Upsilon_{BX'}\pi^{X'}] = S_b + \pi_{B'}\Upsilon_{BX'}\pi^{X'},$$
where $\Upsilon_b = \nabla_b\ln\chi$. Now, (1.3) is equivalent to $S_b = \pi_{A'}\nabla_b\pi^{A'} =: \omega_B\pi_{B'}$, for some spinor $\omega_B$ (see Law 2008, (6.2.6) et seq., for the significance of $S_b$ and $\omega_B$). Thus, one seeks a function $f:U_p \to {\bf R}$ satisfying
$$\pi^{B'}\nabla_{BB'}f = \omega_B,\eqno(3.27.1)$$
$f := \ln\chi^{-1}$. (Notice that (3.27.1) does reduce to (3.1) when $(M,h,Z_{[\pi]})$ is in fact Walker.) This equation has necessary and sufficient integrability condition for local solvability:
$$\pi^{A'}\pi^{B'}\nabla^A_{B'}\omega_A = \omega_A\pi^{B'}\nabla^A_{B'}\pi^{A'}\eqno(3.27.2)$$
(see Law 2008, (6.2.16); Penrose \& Rindler 1986, (7.3.20)). Equation (1.3) is also equivalent to $\pi^{B'}\nabla_{BB'}\pi^{A'} = \eta_B\pi^{A'}$, for some spinor $\eta_B$ (again, see Law 2008, (6.2.6) et seq.), whence (3.27.2) becomes
$$\pi^{A'}\pi^{B'}\nabla^A_{B'}\omega_A = \eta^B\omega_B\pi^{A'}.\eqno(3.27.3)$$
Now, on the one hand
$$\pi^{B'}\nabla^A_{B'}S_a = \pi^{B'}\nabla^A_{B'}(\omega_A\pi_{A'}) =  \pi_{A'}\pi^{B'}\nabla^A_{B'}\omega_A + \omega_A\pi^{B'}\nabla^A_{B'}\pi^{A'} = \pi_{A'}\pi^{B'}\nabla^A_{B'}\omega_A + \eta^A\omega_A\pi_{A'},$$
while on the other,
$$\eqalign{\pi^{B'}\nabla^A_{B'}S_a &= \pi^{B'}\nabla^A_{B'}(\pi_{C'}\nabla_a\pi^{C'})\cr
&= (\pi^{B'}\nabla^A_{B'}\pi_{C'})(\nabla_a\pi^{C'}) + \pi_{C'}\pi^{B'}\nabla^A_{B'}\nabla_a\pi^{C'}\cr
&= \eta^A\pi_{C'}\nabla_a\pi^{C'} + \pi_{C'}\pi^{B'}\nabla^A_{B'}\nabla_{AA'}\pi^{C'}\cr
&= \eta^A\omega_A\pi_{A'} + \pi_{C'}\pi^{B'}\nabla^A_{B'}\nabla_{AA'}\pi^{C'},\cr}$$
whence
$$\pi_{A'}\pi^{B'}\nabla^A_{B'}\omega_A = \pi_{C'}\pi^{B'}\nabla^A_{B'}\nabla_{AA'}\pi^{C'}.\eqno(3.27.4)$$
By Law (2008), (6.2.13)(c),
$$\pi_{C'}\pi^{B'}\nabla^A_{B'}\nabla_{AA'}\pi^{C'} = \eta^B\omega_B\pi_{A'} + 2\tilde\Psi_{A'B'C'D'}\pi^{B'}\pi^{C'}\pi^{D'}.$$
Hence, (3.27.4) becomes
$$\pi_{A'}\pi^{B'}\nabla^A_{B'}\omega_A = \eta^B\omega_B\pi_{A'} + 2\tilde\Psi_{A'B'C'D'}\pi^{B'}\pi^{C'}\pi^{D'},$$
i.e., the integrability condition (3.27.3) is satisfied iff $\tilde\Psi_{A'B'C'D'}\pi^{B'}\pi^{C'}\pi^{D'} = 0$. Thus, iff $[\pi^{A'}]$ is a multiple WPS can one solve (3.27.1) for $f$ on some neighbourhood $U_p$ of $p$. Taking $\Omega = \chi^{-1} = \exp(f)$, then $h = \Omega^2g$ on $U_p$, with $(U_p,g,Z_{[\pi]})$ Walker.\bull
\vskip 24pt
\noindent {\bf 3.28 Remark}\hfil\break
As a corollary of the computations in (3.27), one deduces that for any AS$\alpha$-geometry (real or complex as the relevant computations are valid in both cases), and any LSR $\pi^{A'}$ of $[\pi^{A'}]$,
$$\pi_{A'}\pi^{C'}\Square\pi_{C'} = 2\left[\eta^B\omega_B\pi_{A'} + \tilde\Psi_{A'B'C'D'}\pi^{B'}\pi^{C'}\pi^{D'}\right].$$
In a Walker geometry, each of the terms in this equation vanishes.

Proof. One has, using standard results on curvature (see, for example, Law \& Matsushita 2008, Appendix 2, or Law 2008, Appendix)
$$\eqalign{\pi_{C'}\pi^{B'}\nabla^A_{B'}\nabla_{AA'}\pi^{C'} &= \pi^{C'}\pi^{B'}\nabla_{BB'}\nabla^B_{A'}\pi_{C'}\cr
&= \pi^{C'}\pi^{B'}\left[\Square{}_{B'A'}\pi_{C'} + {1 \over 2}\epsilon_{B'A'}\Square\pi_{C'}\right]\cr
&= \tilde\Psi_{A'B'C'D'}\pi^{B'}\pi^{C'}\pi^{D'} + {1 \over 2}\pi_{A'}\pi^{C'}\Square\pi_{C'}.\cr}$$
Equating with the alternate expression given in the proof of (3.27) yields the assertion.\bull
\vskip 24pt
Thus, the local geometry of any real AS$\alpha$-geometry $(M,h,[\pi^{A'}])$ can be described in a suitable open set $U \subseteq M$, as $(U,\Omega^2g,[\pi^{A'}])$, for some Walker geometry $(U,g,Z_{[\pi]})$. It is straightforward to combine the results of Law \& Matsushita (2008) and Law (2008) on Walker geometry with the conformal rescaling formulae of \S 2 to obtain this description. We will therefore denote the connection and curvature quantities of $(M,h)$ by hatted symbols; unhatted symbols will refer to a (locally) conformally related Walker geometry.

The Weyl curvature spinors $\hat\Psi_{ABCD}$ and $\hat{\tilde\Psi}_{A'B'C'D'}$ of $(U,h,[\pi^{A'}])$ of course coincide with those, $\Psi_{ABCD}$ and $\tilde\Psi_{A'B'C'D'}$ respectively, of $(U,g,Z_{[\pi]})$. The Weyl curvature endomorphisms of $(U,h,[\pi^{A'}])$ are simply scalar multiples of those for $(U,g,Z_{[\pi]})$, $\hat\Psi^{AB}{}_{CD} = \Omega^{-2}\Psi^{AB}{}_{CD}$ and $\hat{\tilde\Psi}^{A'B'}{}_{C'D'} = \Omega^{-2}\tilde\Psi^{A'B'}{}_{C'D'}$, so, for example, the analysis of Law \& Matsushita (2008) between (2.21) and (2.22) applies equally well to $\hat{\tilde\Psi}^{A'B'}{}_{C'D'}$ with exactly the same results for the algebraic classification as for $\tilde\Psi^{A'B'}{}_{C'D'}$ recorded in Law \& Matsushita (2008) 2.6 (the eigenvalues of the Weyl curvature endomorphisms for $(U,h,[\pi^{A'}])$ of course being $\Omega^{-2}$ times those of the corresponding Weyl curvature endomorphism for $(U,g,Z_{[\pi]})$).

For conformal oriented Walker coordinates, when convenient to employ a null tetrad $\{L^a,N^a,M^a,\tilde M^a\}$ and spin frames $\{\hat\alpha^A,\hat\beta^A\}$ and $\{\hat\pi^{A'},\hat\xi^{A'}\}$ for $(U,h,[\pi^{A'}])$, we will exploit those of (3.15--19), in which context $\pi^{A'}$ is always the LSR (unique up to sign) which satisfies (3.5) and which is the first element of the associated Walker primed spin frame for $(U,g,Z_{[\pi]})$. The spin coefficients for $(U,h,[\pi^{A'}])$ with respect to the chosen spin frames (3.15--16) are thus given by (3.21--22). As previously noted, one can substitute for the nonzero spin coefficients of $(U,g,Z_{[\pi]})$ expressions in terms of $a$, $b$, $c$ (of (3.2)) from Law (2008) (5.6). The components of the Weyl curvature spinors of $(U,h,[\pi^{A'}])$ with respect to the chosen spin frames are easily obtained in terms of the components of the Weyl curvature spinors of $(U,g,Z_{[\pi]})$ with respect to the Walker spin frames; for the latter see Law \& Matsushita (2008) (2.20) and (2.25).

At this point, it is convenient to introduce the following notation, defined in terms of oriented Walker coordinates $(u,v,x,y)$ and their associated Walker spin frames $\{\alpha^A,\beta^A\}$ and $\{\pi^{A'},\xi^{A'}\}$:
$$\eqalignno{\delta_A &:= \pi^{A'}\nabla_{AA'} = \alpha_A\tilde m^b\nabla_b - \beta_A\ell^b\nabla_b = \alpha_A\triangle - \beta_AD\cr
&= \alpha_A\partial_v - \beta_A\partial_u,\qquad\hbox{when acting on functions}.&(3.29)\cr}$$
This spinor operator in effect represents the (flat) induced connection, with respect to the Walker geometry, within $\alpha$-surfaces; see Law (2008) (6.2.54) for properties of $\delta_A$.

Turning to the Ricci curvature, a natural geometric condition on $(M,h,[\pi^{A'}])$ is that $[\pi^{A'}]$ be a RPS, i.e., $\hat\Phi_{ABA'B'}\eta^{A'}\eta^{B'} = 0$, for any LSR $\eta^{A'}$ of $[\pi^{A'}]$. (Note that, by Law 2008 (6.2.18), any solution $[\pi^{A'}]$ of (1.3) which is also a RPS is automatically a multiple WPS.) Working in $(U,h)$ with conformal oriented Walker coordinates, since the $\hat\pi^{A'}$ of (3.15) {\sl is} an LSR of $[\pi^{A'}]$, the condition is $\hat\Phi_{00} = \hat\Phi_{10} = \hat\Phi_{20} = 0$. From (2.8), and noting that $[\pi^{A'}]$ is a PS of $\Phi_{ABA'B'}$ (Law \& Matsushita 2008, 2.5),
$$\eqalignno{\hat\Phi_{ABA'B'}\hat\pi^{A'}\hat\pi^{B'} &= \Omega\left[\Upsilon_{AA'}\pi^{A'}\Upsilon_{BB'}\pi^{B'} - \pi^{B'}\pi^{A'}\nabla_{AA'}\Upsilon_{BB'}\right]\cr
&= \Omega\left[(\delta_A\omega)(\delta_B\omega) - \pi^{B'}\delta_A\Upsilon_{BB'}\right]&(3.30)\cr
&= \Omega\left[(\delta_A\omega)(\delta_B\omega) - \delta_A\delta_B\omega\right],\cr}$$
since, by Law (2008) (5.8), the Walker spin frames are parallel on $\alpha$-surfaces (in particular $\delta_A\pi^{B'} = 0$).
\vskip 24pt
\noindent {\bf 3.31 Lemma}\hfil\break
For a real AS$\alpha$-geometry $(M,h,[\pi^{A'}])$, $[\pi^{A'}]$ is a RPS iff, with respect to any conformal Walker coordinates $(u,v,x,y)$ on $U \subseteq M$ for $(M,h,[\pi^{A'}])$ (i.e., $(U,h) = (U,\Omega^2g)$ and $(u,v,x,y)$ are Walker coordinates on $U$ for the Walker geometry $(U,g,Z_{[\pi]})$), $\Omega^{-1}$ is affine as a function of $u$ and $v$, i.e.,
$$\Omega(u,v,x,y) = \bigl[M(x,y)u + N(x,y)v + K(x,y)\bigr]^{-1},$$
for some functions $M$, $N$, and $K$ of $(x,y)$.

Proof. Without loss of generality, one can suppose the conformal Walker coordinates are oriented. From (3.30), $[\pi^{A'}]$ is a PS of $\hat\Phi_{ABA'B'}$ iff 
$$\delta_A\delta_B\omega = (\delta_A\omega)(\delta_B\omega).\eqno(3.31.1)$$
But (3.31.1) is equivalent to
$$\Omega\delta_A\delta_B\Omega = 2(\delta_A\Omega)(\delta_B\Omega),\eqno(3.31.2)$$
which is equivalent to
$$\delta_A\delta_B\left(\Omega^{-1}\right) = 0.\eqno(3.31.3)$$
By Law (2008) (6.2.54), (3.31.3) is equivalent to the assertion of the lemma.
One can also obtain the desired result from the spin coefficient field equations Law (2008) (3.4). For example, Law (2008) (3.4a), in conjunction with (3.21--22), yields
$$-\hat D\hat\rho = \hat\rho^2 - \hat\rho(\hat\epsilon+\hat{\tilde\epsilon}) + \hat\Phi_{00}.$$
Substituting for these nonzero spin coefficients the expressions in (3.21--22), and noting $\hat D = D$ on functions for (3.15--19), results in $-D^2\omega = (D\omega)^2 - (D\omega)(2D\omega) + \hat\Phi_{00}$, i.e.,
$$\hat\Phi_{00} = -D^2\omega + (D\omega)^2,$$
which is the result obtained upon transvecting (3.30) by $\hat\alpha^A\hat\alpha^B$. Similarly, the equations resulting from transvecting (3.30) by $\hat\alpha^A\hat\beta^B$ (or $\hat\beta^A\hat\alpha^A$) and $\hat\beta^A\hat\beta^B$ also result, respectively, from Law (2008) (3.4) (\~ c) (or (\~d)) and (\~e).\bull
\vskip 24pt
Note that one can always introduce Frobenius coordinates $(p,q,x,y)$ for the integrable $\alpha$-distribution $Z_{[\pi]}$ in $(M,h,[\pi^{A'}])$, with $x$ and $y$ constant on $\alpha$-surfaces. With respect to such Frobenius coordinates, the metric $h$ takes the form
$$\left(h_{\rm ab}\right) = \pmatrix{{\bf 0}_2&V\cr {^\tau\! V}&W\cr}.$$
for some $V$, $W \in {\bf R}(2)$, with $W$ symmetric. Since $dx$ and $dy$ vanish on $\alpha$-surfaces, for any LSR $\pi^{A'}$ one can write
$$dx = \mu_A\pi_{A'} \hskip 1in dy = \nu_A\pi_{A'} \hskip .5in \nu_A\mu_A \not= 0.$$
The vector fields $U^a := \mu^A\pi^{A'}$ and $V^a := \nu^A\pi^{A'}$ span $Z_{[\pi]}$. Repeating the computation in Law \& Matsushita (2008) 2.3 yields $U^b\nabla_bV^a = (\nu \cdot \mu)S^a = -V^b\nabla_bU^a$, whence $[U,V]^a = 2(\nu \cdot \mu)S^a$, explicitly revealing how $S^a$ blocks the construction of Walker coordinates. If desired, as in Law \& Matsushita (2008) 2.4, there is a unique, up to sign, LSR $\pi^{A'}$ such that $\nu \cdot \mu = \pm 1$. 

If one wishes to retain the freedom to choose any LSR, however, one can choose spin frames so that $o^{A'} = \pi^{A'}$, $o^A \propto \mu^A$ and $\iota^A \propto \nu^A$ (assuming the coordinates $x$ and $y$ are oriented so that $\nu \cdot \mu > 0$), whence $\ell_a \propto dx$, $\tilde m_a \propto dy$. Doing so, write
$$\ell_a = Ldx \hskip 1.25in \tilde m_a = Mdy.$$
Now, completing $\ell_a$ and $\tilde m_a$ to a null tetrad $\{\ell_a,n_a,m_a,\tilde m_a\}$,
$$m_a = A_1dp + B_1dq + C_1dx + D_1dy \hskip 1.25in n_a = A_2dp + B_2dq + C_2dx + D_2dy,$$
where the $A_i$, $B_i$, $C_i$ and $D_i$, $i=1$, 2, are functions of $(p,q,x,y)$. With $x$ and $y$ held constant, each of $\phi_1 := A_1dp+B_1dq$ and $\phi_2 := A_2dp + B_2dq$ are one-forms on the $\alpha$-surface in question. They are linearly dependent iff $\left\vert{A_1 \atop A_2}{B_1 \atop B_2}\right\vert = 0$, in which case $B_2m_a - B_1n_a \in \langle dx, dy \rangle_{\bf R} = \langle\ell_a, \tilde m_a\rangle_{\bf R}$. Linear independence of the null tetrad one-forms then entails $B_1 = B_2 = 0$. Similarly, $A_1 = A_2 = 0$. But then $m_a$ and $n_a$ would be linear combinations of $\ell_a$ and $\tilde m_a$. Hence, $\phi_1$ and $\phi_2$ are linearly independent. Now any one-dimensional distribution is integrable. In a two-dimensional surface, such a distribution is also of co-dimension one. By the differential form version of Frobenius' theorem, any one-form on a two-surface is therefore proportional to a gradient. Hence, one can write
$$\phi_1 = A_1dp + B_1dq =: \zeta dv \hskip 1.25in \phi_2 = A_2dp + B_2dq =: \xi du,$$
where $\zeta$, $\xi$, $u$ and $v$ are functions of $(p,q,x,y)$. Now, on any $\alpha$-surface,
$$du = u_pdp + u_qdq = {A_2 \over \xi}dp + {B_2 \over \xi}dq \hskip 1.25in dv = v_pdp + v_qdq = {A_1 \over \zeta}dp + {B_1 \over \zeta}dq,$$
whence the Jacobian of $(u,v,x,y)$ as a function of $(p,q,x,y)$ is nonsingular. Hence, one can use $(u,v,x,y)$ as Frobenius coordinates. Without restriction to a given $\alpha$-surface, one has
$$\zeta dv = A_1dp + B_1dq + C_2dx + D_3dy \hskip 1.25in \xi du = A_2dp + B_2dq +C_4dx + D_4dy$$
whence, for some functions $\hat C_1$, $\hat C_2$, $\hat D_1$, $\hat D_2$, $A$, $B$, $C$, $D$, 
$$m_a = \zeta dv + \hat C_1dx + \hat D_1dy = \zeta(dv + Cdx + Ddy) \hskip .5in n_a = \xi du + \hat C_2dx + \hat D_2dy = \xi(du + Adx + Bdy)$$
and
$$\eqalign{h_{ab} &= 2[\ell_{(a}n_{b)} - m_{(a}\tilde m_{b)}]\cr
&= 2\bigl[L\xi dudx - M\zeta dvdy + L\xi dx(Adx + Bdy) - M\zeta dy(Cdx+Ddy)\bigr]\cr
&= 2\bigl[\phi dudx + \psi dvdy + \phi dx(Adx + Bdy) + \psi dy(Cdx + Ddy)\bigr]\cr}$$
where $\phi := L\xi$ and $\psi := M\zeta$, i.e.,
$$h_{\bf ab} = \pmatrix{0&0&\phi&0\cr 0&0&0&\psi\cr \phi&0&2\phi A&\phi B + \psi C\cr 0&\psi&\phi B + \psi C&2\psi D\cr}.\eqno(3.32)$$
\vskip 6pt
(3.32) is a generalization of Walker's form of the metric when the integrable $\alpha$-distribution is not assumed to be parallel.  Pleba\'nski and Robinson (1976), (1977) derived this coordinate form in the context of complex general relativity to study algebraically special Weyl curvature spinors. Their approach was developed into the {\sl hyperheavenly formalism} for such complex space-times, in which the Einstein condition is reduced to a single PDE, the {\sl hyperheavenly equation}, see Finley et al. (1976), Boyer et al. (1980). The hyperheavenly formalism is obviously equally applicable to real (four-dimensional) neutral geometry. The derivation of the hyperheavenly equation involves making further judicious changes of coordinates (in particular, the assumption of algebraic degeneracy of the relevant Weyl curvature spinor permits a change of coordinates which results in $\phi = \psi$) and other choices which simplify the expression of the Einstein condition. Two separate cases arise: the {\sl expanding} and {\sl nonexpanding} cases. In the real neutral context, the nonexpanding case corresponds to Walker geometry; an independent derivation of the hyperheavenly equation in this context, utilizing the spinor approach of Law \& Matsushita (2008), was presented in Law (2008), (6.2.45--63). Real AS$\alpha$-geometries which are not Walker correspond to the expanding case in the hyperheavenly formalism. The approach to characterizing the local geometry of real AS$\alpha$-nonWalker geometries developed in this paper, based on the local conformally Walker structure, takes full advantage of Walker geometry. Though our approach has been independent of the hyperheavenly formalism, further development closely parallels the derivation of the hyperheavenly equation.

We next show it is possible to choose conformal (oriented) Walker coordinates which simplify the result of (3.31).
\vskip 24pt
\noindent {\bf  3.33 Lemma}\hfil\break
For a real AS$\alpha$-geometry $(M,h,[\pi^{A'}])$ for which $[\pi^{A'}]$ is a RPS, every $p \in M$ has a neighbourhood $U$ admitting conformal (oriented) Walker coordinates $(U,V,X,Y)$ (i.e., $(U,h) = (U,\Omega^2g)$ for which $(U,V,X,Y)$ are (oriented) Walker coordinates for the Walker geometry $(U,g,Z_{[\pi]})$) such that 
$$\Omega(U,V,X,Y) = (MU+NV)^{-1}$$
for constants $M$ and $N$.

Proof. By (3.31), one can choose conformal (oriented) Walker coordinates $(u,v,x,y)$ for $(M,h,[\pi^{A'}])$ on a neighbourhood of $p$ so that $\Omega(u,v,x,y) = M(x,y)u + N(x,y)v + K(x,y)$. We seek new coordinates $(U,V,X,Y)$, on a possibly smaller neighbourhood, which are still conformal (oriented) Walker coordinates but for which the desired result holds. Following Pleba\'nski \& Robinson (1976), consider the following coordinate transformation:
$$\vcenter{\openup1\jot \halign{$\hfil#$&&${}#\hfil$&\qquad$\hfil#$\cr
U &= [uY_y - vY_x + P(x,y)]H(x,y) & X &= X(x,y)\cr
V &= [vX_x - uX_y + Q(x,y)]H(x,y) & Y &= Y(x,y),\cr}}$$
\vskip 6pt
where $X$, $Y$, $P$, $Q$, $H$ are yet-to-be-determined functions. The Jacobian for this transformation is
$$J := {\partial(U,V,X,Y) \over \partial(u,v,x,y)} = \pmatrix{D&G\cr {\bf 0}_2&F\cr},$$
where
$$D = H\pmatrix{Y_y&-Y_x\cr -X_y&X_x\cr} \hskip 1.25in F = \pmatrix{X_x&X_y\cr Y_x&Y_y\cr}.$$
As $(u,v,x,y)$ are Walker coordinates for $(U,g,Z_{[\pi]})$, so also are $(U,V,X,Y)$ iff $F = {^\tau\! D}^{-1}$, see Law \& Matsushita (2008) (1.2), which is the case iff $H(x,y) = (X_xY_y - Y_xX_y)^{-1}$. This condition therefore determines $H$ in terms of $X$ and $Y$.

For constants $K_1$ and $K_2$,
$$K_1U + K_2V = (K_1Y_y - K_2X_y)Hu + (K_2X_x - K_1Y_x)Hv + (K_1P + K_2Q)H.$$
To prove the lemma, one must therefore show that one can choose $X$, $Y$, $P$ and $Q$ so that $(X_xY_y - Y_xX_y) \not= 0$ and
$$\pmatrix{M(x,y)\cr N(x,y)\cr} = H\pmatrix{Y_y&-X_y\cr -Y_x&X_x\cr}.\pmatrix{K_1\cr K_2\cr} \hskip 1in H(K_1P + K_2Q) = K.$$
We may construe this problem as that of finding a change of coordinates $\bigl(X(x,y),Y(x,y)\bigr)$ on some neighbourhood of the origin in ${\bf R}^2$ so that the given vector field $V(x,y) := M(x,y)\partial_x + N(x,y)\partial_y$ transforms to the constant vector field $Z(x,y) := K_1\partial_X + K_2\partial_Y$ and so that $K_1P+K_2Q = KH^{-1}$. Given that $V(x,y)$ is nonzero at the origin, it is always possible to find such a change of coordinates; once $X$ and $Y$ are chosen and $H$ thereby determined, the final condition has many solutions for $P$ and $Q$.\bull
\vskip 24pt
The value of (3.33) is that such conformal (oriented) Walker coordinates simplify subsequent analysis of the Ricci curvature. This choice corresponds to a similar manoeuvre in the hyperheavenly formalism, e.g., in Finley et al. (1976), between (3.1) and (3.3), and in Boyer et al. (1980) at (5.21). In the remainder of this section, we will therefore employ such conformal oriented Walker coordinates $(U,V,X,Y)$ and will write 
$$\Omega = (MU+NV)^{-1} \hskip 1.25in {\cal W} = \pmatrix{A&C\cr C&B\cr},\eqno(3.34)$$
where $M$ and $N$ are constants and $\cal W$ is the analogue of $W$ in (3.2--3) for $(U,V,X,Y)$. It is convenient to record:
$${\partial^{p+q}\Omega \over \partial U^p\partial V^q} = (-1)^{p+q}(p+q)!M^pN^q\Omega^{p+q+1}; \qquad \delta_A\Omega = -\Omega^2(\alpha_AN - \beta_AM); \qquad \delta_A\delta_B\Omega = 2\Omega^{-1}(\delta_A\Omega)(\delta_B\Omega_B).\eqno(3.35)$$

As an example of the utility of (3.33), consider the Ricci scalar curvature $S = -24\Lambda$. From (2.8), one needs to evaluate $\Square\Omega$ in the Walker geometry $(U,g,Z_{[\pi]})$. From Law \& Matsushita (2008) (3.9), for any function $F$, with respect to Walker coordinates $(U,V,X,Y)$ and with notation as in (3.34),
$$\Square F = -AF_{UU} - 2CF_{UV} - BF_{VV} + 2F_{XU} + 2F_{YV} - (A_U + C_V)F_U - (B_V + C_U)F_V.$$
In particular,
$$\Square\Omega = -2(AM^2 + 2CMN + BN^2)\Omega^3 + (A_U + C_V)M\Omega^2 + (B_V + C_U)N\Omega^2,$$
whence, noting from Law \& Matsushita (2008) A1.6 that $-24\Lambda = S = A_{UU} + B_{VV} + 2C_{UV}$,
$$\eqalignno{\hat\Lambda &= \Omega^{-2}\left[\Lambda - {1 \over 2}(AM^2 + 2CMN + BN^2)\Omega^2 + {1 \over 4}\bigl((A_U + C_V)M + (B_V + C_U)N\bigr)\Omega\right]\cr
&= -{\Omega^{-5} \over 24}\left[\bigl(A_{UU}\Omega^3 - 6A_UM\Omega^4 + 12AM^2\Omega^5\bigr) + \bigl(B_{VV}\Omega^3 - 6B_VN\Omega^4 + 12BN^2\Omega^5\bigr)\right]\cr
&\qquad -{\Omega^{-5} \over 24}\left[2\bigl(C_{UV}\Omega^3 - 3(C_VM+C_UN)\Omega^4 + 12CMN\Omega^5\bigr)\right]\cr
&= -{\Omega^{-5} \over 24}\left[\bigl(A\Omega^3\bigr)_{UU} + \bigl(B\Omega^3\bigr)_{VV} + 2\bigl(C\Omega^3\bigr)_{UV}\right]\cr
&= -{\Omega^{-5} \over 24}\delta_{\bf A}\delta_{\bf B}H^{\bf AB}\qquad\hbox{where}\qquad(H^{\bf AB}) := \Omega^3{\cal W}.&(3.36)\cr}$$
As a result, it is convenient to label $\cal W$ with concrete superscript `spinor' indices.

Granted that $[\pi^{A'}]$ is a PS of $\hat\Phi_{ABA'B'}$, it is natural to investigate the condition for $[\pi^{A'}]$ to be a multiple RPS; namely $\hat\Phi_{ABA'B'}\eta^{B'} = 0$, for any LSR $\eta^{A'}$ of $[\pi^{A'}]$. Once again, work on an open set $U$ where $(U,h) = (U,\Omega^2g)$, $(U,g,Z_{[\pi]})$ is Walker, and on which $(u,v,x,y)$ are oriented Walker coordinates for $(U,g,Z_{[\pi]})$. We continue to employ for spin frames for $(U,h)$ those of (3.15--19) constructed from the Walker spin frames $\{\alpha^A,\beta^A\}$ and $\{\pi^{A'},\xi^{A'}\}$ associated to the oriented Walker coordinates $(u,v,x,y)$ for $(U,g,Z_{[\pi]})$. Hence, the condition for $[\pi^{A'}]$ to be a multiple PS of $\hat\Phi_{ABA'B'}$ may be expressed as $\hat\Phi_{ABA'B'}\hat\pi^{B'} = 0$, i.e., in addition to $\hat\Phi_{00} = \hat\Phi_{10} = \hat\Phi_{20} = 0$, one also requires $\hat\Phi_{01} = \hat\Phi_{11} = \hat\Phi_{21} = 0$. These quantities may be computed from Law (2008), (3.4) \& (5.6), together with (3.21--22); but it is more efficient to exploit (2.8).
\vskip 24pt
\noindent {\bf 3.37 Lemma}\hfil\break 
With assumptions and notation as in the previous paragraph, with $B_{AB} := -2\Phi_{ABA'B'}\xi^{A'}\pi^{B'}$, see Law \& Matsushita (2008) (2.32--33), and $\varsigma_A$ as defined in (A6), then
$$\hat\Phi_{AB1'0'} := \hat\Phi_{ABA'B'}\hat\xi^{A'}\hat\pi^{B'} = \Omega^{-1}\left[-{1 \over 2}B_{AB} -\Omega\delta_{(A}[\Omega^{-2}\varsigma_{B)}\Omega]\right].\eqno(3.37.1)$$

Proof. Recalling, from Law (2008) (5.8), that the Walker spin frames are parallel with respect to $\delta_A$, one computes from (2.8)
$$\eqalign{\hat\Phi_{AB1'0'} &= \Omega^{-1}\left[\Phi_{ABA'B'} + \Upsilon_{A(A'}\Upsilon_{B')B} - \nabla_{A(A'}\Upsilon_{B')B}\right]\xi^{A'}\pi^{B'}\cr
&= {\Omega^{-1} \over 2}\left[-B_{AB} + 2[\varsigma_{(A}\omega][\delta_{B)}\omega] - \pi^{B'}\varsigma_A\Upsilon_{BB'} - \delta_A\varsigma_B\omega\right]\cr
&= \Omega^{-1}\left[-{1 \over 2}B_{AB} + 2\Omega^{-2}[\varsigma_{(A}\Omega][\delta_{B)}\Omega] - {\Omega^{-1} \over 2}(\pi^{B'}\varsigma_A\nabla_{BB'}\Omega + \delta_A\varsigma_B\Omega)\right]\cr
&= \Omega^{-1}\left[-{1 \over 2}B_{AB} + 2\Omega^{-2}[\varsigma_{(A}\Omega][\delta_{B)}\Omega] - {\Omega^{-1} \over 2}\bigl(\varsigma_A\delta_B\Omega - (\varsigma_A\pi^{B'})\nabla_{BB'}\Omega + \delta_A\varsigma_B\Omega\bigr)\right].\cr}$$
The term in round brackets does not appear to be symmetric in $A$ and $B$ but in fact it is (as can be seen by computing the unsymmetrized version of (A7)), as of course it must be. Symmetrizing over $A$ and $B$ in the previous equation and using (A7) yields
$$\hat\Phi_{AB1'0'} = \Omega^{-1}\left[-{1 \over 2}B_{AB} + 2\Omega^{-2}[\varsigma_{(A}\Omega][\delta_{B)}\Omega] - \Omega^{-1}\delta_{(A}\varsigma_{B)}\Omega\right],$$
which proves the assertion.\bull
\vskip 24pt
Note that as the Walker spin frames are parallel with respect to $\delta_A$, the components of $\hat\Phi_{AB1'0'}$ with respect to the spin frames (3.15--19) are just
$$\hat\Phi_{{\bf AB}1'0'} = \Omega^{-2}\left[-{1 \over 2}B_{\bf AB} -\Omega\delta_{({\bf A}}\Omega^{-2}\varsigma_{\bf B)}\Omega\right],\eqno(3.38)$$
where the components on the right-hand side of (3.38) refer to the Walker spin frames.

From Law \& Matsushita (2008) (2.32--33) \& A1.8, one can write $B_{AB}$ in the form
$$B_{AB} = {1 \over 2}\delta_{(A}\delta_{\vert{\bf C}\vert}W^{\bf BC}\epsilon_{B)D}\,\epsilon_{\bf B}{}^D = - {1 \over 2}\delta_{(A}\delta_{\vert{\bf C}\vert}W^{\bf BC}\epsilon_{B)}{}^{\bf D}\epsilon_{\bf BD}.\eqno(3.39)$$
\vskip 24pt
\noindent {\bf 3.40 Lemma}\hfil\break
With assumptions as in (3.37), but utilizing the conformal oriented Walker coordinates $(U,V,X,Y)$ of (3.33) rather than arbitrary conformal oriented Walker coordinates,
$$\hat\Phi_{AB1'0'} = {\Omega^{-3} \over 4}\delta_{\bf C}\left[\Omega^2\delta_{(A}{\cal W}^{\bf KC}\epsilon_{B)}{}^{\bf D}\right]\epsilon_{\bf KD}.$$

Proof. From Law \& Matsushita (2008) (2.11), 
$$\varsigma_A  = \alpha_A\left[\partial_X  -  {1 \over  2}(AD+C\triangle)\right]  -  \beta_A\left[{1 \over 2}(CD+B\triangle)-\partial_Y\right] = \alpha_A\partial_X + \beta_A\partial_Y  + {1  \over 2}\epsilon_A{}^{\bf F}\epsilon_{\bf  FC}{\cal W}^{\bf  CD}\delta_{\bf D}.\eqno(3.40.1)$$
With respect to $(U,V,X,Y)$, $\Omega$ is as in (3.34), whence independent of $X$ and $Y$, so
$$\varsigma_B\Omega = {1 \over 2}\epsilon_{\bf FC}{\cal W}^{\bf CD}\delta_{\bf D}\Omega\epsilon_B{}^{\bf F}.$$
Hence,
$$\eqalign{\delta_A(\Omega^{-2}\varsigma_B\Omega) &= -{1 \over 2}\left[\delta_A(\Omega^{-2}{\cal W}^{\bf CD}\delta_D\Omega)\right]\epsilon_{\bf CF}\epsilon_B{}^{\bf F}\epsilon_{\bf D}{}^D\cr
&= -{1 \over 2}\left[-2\Omega^{-3}{\cal W}^{\bf CD}(\delta_D\Omega)(\delta_A\Omega) + \Omega^{-2}{\cal W}^{\bf CD}\delta_A\delta_D\Omega + \Omega^{-2}(\delta_A{\cal W}^{\bf CD})(\delta_D\Omega)\right]\epsilon_{\bf CF}\epsilon_B{}^{\bf F}\epsilon_{\bf D}{}^D\cr
&= -{1 \over 2}\left[\Omega^{-2}(\delta_A{\cal W}^{\bf CD})(\delta_D\Omega)\right]\epsilon_{\bf CF}\epsilon_B{}^{\bf F}\epsilon_{\bf D}{}^D,\qquad\hbox{upon using (3.35)}.\cr}$$
Substituting this last expression and (3.39) into (3.37.1) yields,
$$\eqalign{\hat\Phi_{AB1'0'} &= {\Omega^{-1} \over 4}\bigl[\delta_{(A}\delta_{\vert{\bf C}\vert}{\cal W}^{\bf KC}\epsilon_{B)}{}^{\bf D} + 2\Omega^{-1}[\delta_{(A}{\cal W}^{\bf KC}]\epsilon_{B)}{}^{\bf D}\delta_{\bf C}\Omega\bigr]\epsilon_{\bf KD}\cr
&= {\Omega^{-3} \over 4}\bigl[\Omega^2\delta_{\bf C}\delta_{(A}{\cal W}^{\bf KC}\epsilon_{B)}{}^{\bf D} + 2\Omega[\delta_{(A}{\cal W}^{\bf KC}\epsilon_{B)}{}^{\bf D}]\delta_{\bf C}\Omega\bigr]\epsilon_{\bf KD}\cr
&= {\Omega^{-3} \over 4}\delta_{\bf C}\left[\Omega^2\delta_{(A}{\cal W}^{\bf KC}\epsilon_{B)}{}^{\bf D}\right]\epsilon_{\bf KD}.\bull\cr}$$
\vskip 24pt
Hence, with respect to conformal oriented Walker coordinates $(U,V,X,Y)$ as in (3.33), (3.38) specializes to
$$\hat\Phi_{{\bf AB}1'0'} = -{\Omega^{-4} \over 4}\delta_{\bf C}(\Omega^2\delta_{({\bf A}}{\cal W}^{\bf KC}\epsilon_{\bf B)K}).\eqno(3.41)$$
Specifically:
$$\displaylines{\hat\Phi_{001'0'} = -{\Omega^{-2} \over 4}\left[C_{UU} + B_{UV} - 2\Omega(MC_U + NB_U)\right] \qquad \hat\Phi_{111'0'} = {\Omega^{-2} \over 4}\left[A_{UV} + C_{VV} - 2\Omega(MA_V + NC_V)\right]\cr
\hat\Phi_{011'0'} = -{\Omega^{-2} \over 8}\left[B_{VV} - A_{UU} + 2\Omega\bigl((A_U-C_V)M + (C_U-B_V)N\bigr)\right].\cr}$$
\vskip 24pt
\noindent{\bf 3.42 Lemma}\hfil\break
With assumptions as in (3.40) and exploiting the notation of (A10), if the Ricci scalar curvature $\hat S$ of the real AS$\alpha$-geometry $(M,h,[\pi^{A'}])$ is constant then, on a suitable neighbourhood,
$${\cal W}^{\bf AB} = \Omega^{-3}\delta^{\bf (A}F^{\bf B)} + {\hat S \over 12\tau^2}K^{\bf A}K^{\bf B},\eqno(3.42.1)$$
for some pair of functions $F^{\bf B}$, and
$$\hat\Phi_{AB1'0'} = {\Omega^{-3} \over 8}\delta_A\delta_B\bigl(\Omega\delta_{\bf C}(\Omega^{-2}F^{\bf C})\bigr).\eqno(3.42.2)$$

Proof. The proof is adapted from Finley \& Pleba\'nski (1976). Notice that (3.36) reads $\delta_{\bf A}\delta_{\bf B}H^{\bf AB} = \hat S\Omega^5$, where $H^{\bf AB} = \Omega^3{\cal W}^{\bf AB}$. With $\hat S$ constant, the linear system of PDEs $\delta_{\bf A}E^{\bf A} = \hat S\Omega^5$ has general solution $E^{\bf A} = \delta^{\bf A}H - (\hat S/4\tau)\Omega^4K^{\bf A}$, for arbitrary functions $H$. Thus, $\delta_{\bf B}H^{\bf AB} = \delta^{\bf A}H - (\hat S/4\tau)\Omega^4K^{\bf A}$, for some specific $H$. Write $H =: (1/2)\delta_{\bf B}Z^{\bf B}$ for some (nonunique) quantities $Z^{\bf B}$. With 
$$F^{\bf AB} := H^{\bf AB} - \delta^{\bf (A}Z^{\bf B)} - {\hat S\Omega^3\over 12\tau^2}K^{\bf A}K^{\bf B},$$
$\delta_{\bf B}F^{\bf AB} = 0$. By the Poincar\'e lemma, $F^{\bf AB} = \delta^{\bf B}G^{\bf A}$, for some functions $G^{\bf A}$. But, as $F^{\bf AB} = F^{\bf BA}$, then $\delta^{\bf A}G^{\bf B} = \delta^{\bf B}G^{\bf A}$, i.e., $-\partial_U(G^0) = \partial_V(G^1)$, whence, by the Poincar\'e lemma again, there exists a function $K$ such that $G^{\bf A} = \delta^{\bf A}K$. So, $F^{\bf AB} = \delta^{\bf A}\delta^{\bf B}K$, for some function $K$. Hence,
$${\cal W}^{\bf AB} = \Omega^{-3}H^{\bf AB} = \Omega^{-3}(\delta^{\bf A}\delta^{\bf B}K + \delta^{\bf (A}Z^{\bf B)}) + {\hat S \over 12\tau^2}K^{\bf A}K^{\bf B} = \Omega^{-3}\delta^{\bf (A}F^{\bf B)} + {\hat S \over 12\tau^2}K^{\bf A}K^{\bf B},$$
where $F^{\bf B} = \delta^{\bf B}K + Z^{\bf B}$. 

In the hyperheavenly formalism, Pleba\'nski and co-workers showed that it is possible to choose coordinates in (3.32) so that $\phi=\psi$, and then, assuming that $[\pi^{A'}]$ is a RPS, so that $\phi$ is independent of the $x$ and $y$ coordinates and linear in the $u$ and $v$ coordinates. In effect, one can take $\phi$ to have the same functional form as $\Omega^{-1}$ has as a function of the conformal oriented Walker coordinates $(U,V,X,Y)$. Formally identifying their $\phi$ with $\Omega^{-1}$, one sees that the expression given for their quantity $C_{12\dot A\dot B}$ (see Finley and Pleba\'nski 1976, p. 2213, and their Appendix B), which is effectively equivalent to $\hat\Phi_{AB1'0'}$, is of the same form, up to a multiplicative factor, as (3.41). Hence, as the second term in (3.42.1) is independent of $U$ and $V$, by the same argument as given in Finley and Pleba\'nski (1976), p. 2213, substitution of (3.42.1) into (3.41) yields (3.42.2), noting (3.15) when taking spin frame components of (3.42.2) (which result we have also confirmed directly).\bull
\vskip 24pt
\noindent {\bf 3.43 Corollary}\hfil\break
For a real AS$\alpha$-geometry $(M,h,[\pi^{A'}])$ for which $[\pi^{A'}]$ is a RPS and the Ricci scalar curvature is constant, $[\pi^{A'}]$ is a multiple RPS iff, with respect to conformal oriented Walker coordinates $(U,V,X,Y)$ satisfying (3.33), there is a function $\vartheta(U,V,X,Y)$ such that
$${\cal W}^{\bf AB} = \Omega^{-3}\delta^{\bf (A}[\Omega^2\delta^{\bf B)}\vartheta] + {2S\Omega^{-2} + \hat S \over 12\tau^2}K^{\bf A}K^{\bf B},\eqno(3.43.1)$$
where $S$ is the Ricci scalar curvature of the Walker geometry $(U,g,Z_{[\pi]})$, with $h = \Omega^2g$ on the neighbourhood $U$, and $S\Omega$ is a function of $(X,Y)$ only.

Proof. Continuing to follow Finley \& Pleba\'nski (1976), (3.42.2) is zero iff $\delta_A\delta_B\bigl(\Omega\delta_{\bf C}(\Omega^{-2}F^{\bf C})\bigr) = 0$, i.e., iff
$$\Omega\delta_{\bf C}(\Omega^{-2}F^{\bf C}) = P(X,Y)U + Q(X,Y)V + R(X,Y),\eqno(3.43.2)$$
for some functions $P$, $Q$, $R$ of $(X,Y)$. Put 
$$(H_{\bf A}) := {1 \over  2}\pmatrix{P\cr Q\cr} \qquad (T^{\bf A}) := \pmatrix{U\cr V\cr}\qquad\hbox{whence}\qquad \delta^{\bf A}T^{\bf B} = \epsilon^{\bf AB}.\eqno(3.43.3)$$
With $\beta^{\bf A}(X,Y)$ functions such that $\beta^{\bf A}J_{\bf A} = -R/2$, define
$$L^{\bf A} := F^{\bf A} - {1 \over 2\tau}K^{\bf B}H_{\bf B}T^{\bf A} + \beta^{\bf A}.\eqno(3.43.4)$$
Since $\delta_{\bf A}\beta^{\bf B} = 0$ and $\delta^{\bf (A}T^{\bf B)} = 0$, then $\delta^{\bf (A}L^{\bf B)} = \delta^{\bf (A}F^{\bf B)}$, whence one can replace $F^{\bf B}$ by $L^{\bf B}$ in (3.42.1). In place of (3.43.2), one obtains
$$\eqalignno{\Omega\delta_{\bf C}(\Omega^{-2}L^{\bf C}) &= \Omega\delta_{\bf C}\left[\Omega^{-2}(F^{\bf C} - {1 \over 2\tau}K^{\bf B}H_{\bf B}T^{\bf C} + \beta^{\bf C})\right]\cr
&= (2T^{\bf B}H_{\bf B} + R) - \Omega\left[{K^{\bf B}H_{\bf B} \over 2\tau}\delta_{\bf C}(\Omega^{-2}T^{\bf C}) - \delta_{\bf C}(\Omega^{-2}\beta^{\bf C})\right]\qquad\hbox{substituting (3.43.2)}\cr
&= 2T^{\bf B}H_{\bf B} + R - \Omega\left[-2\Omega^{-3}({K^{\bf B}H_{\bf B} \over 2\tau}T^{\bf C} - \beta^{\bf C})\delta_{\bf C}\Omega + \Omega^{-2}{K^{\bf  B}H_{\bf B} \over \tau}\right]\cr
&\qquad(\hbox{ as } \delta_{\bf C}\beta^{\bf C} = 0,\ \delta_{\bf C}T^{\bf C} =  \epsilon_{\bf C}{}^{\bf C} = 2)\cr
&= 2T^{\bf B}H_{\bf B} + R - {K^{\bf B}H_{\bf B} \over \tau}T^{\bf C}J_{\bf C} + 2\beta^{\bf C}J_{\bf C} - {K^{\bf B}H_{\bf B}T^{\bf C}J_{\bf C} \over \tau}\qquad(\hbox{noting that }\Omega^{-1} = T^{\bf C}J_{\bf C})\cr
&= 2T^{\bf B}H_{\bf A}\left(\delta^{\bf A}_{\bf B} - {K^{\bf A}J_{\bf B} \over \tau}\right)\cr
&= -{2T^{\bf B}K_{\bf B}J^{\bf A}H_{\bf A} \over \tau},\qquad\hbox{using (A10)},\cr
&= {2\mu T^{\bf B}K_{\bf B} \over \tau},\qquad\hbox{where}\qquad \mu(X,Y) :=  - J^{\bf A}H_{\bf A},&(3.43.5)\cr}$$
which is simpler than (3.43.2). Thus, replacing $F^{\bf C}$ by $L^{\bf C}$ leaves (3.42.1) valid but now $\Phi_{AB1'0'}$ vanishes iff $\Omega\delta_{\bf C}(\Omega^{-2}L^{\bf C}) = 2(\mu/\tau)T^{\bf B}K_{\bf B}$. As this equation is again of the form $\delta_{\bf C}\phi^{\bf C} = \psi$, one need only find a particular solution to write down the general solution, which is
$$L^{\bf C} = \Omega^2\delta^{\bf C}\vartheta + {\mu \over \tau^2}T^{\bf B}K_{\bf B}K^{\bf C},\eqno(3.43.6)$$
where $\vartheta$ is an arbitrary function of $(U,V,X,Y)$. Noting that $K^{\bf (B}\delta^{\bf A)}T^{\bf D}K_{\bf D} = K^{\bf (B}\epsilon^{\bf A)D}K_{\bf D} = K^{\bf A}K^{\bf B}$, substituting (3.43.6) for $F^{\bf C}$ in (3.42.1) yields
$$\eqalignno{{\cal W}^{\bf AB} &= \Omega^{-3}\delta^{\bf (A}[\Omega^2\delta^{\bf B)}\vartheta] + {\mu\Omega^{-3} \over \tau^2}K^{\bf A}K^{\bf B} + {\hat S \over 12\tau^2}K^{\bf A}K^{\bf B}\cr
&= \Omega^{-3}\delta^{\bf (A}[\Omega^2\delta^{\bf B)}\vartheta] + {12\mu\Omega^{-3} + \hat S \over 12\tau^2}K^{\bf A}K^{\bf B}.&(3.43.7)\cr}$$
At this point, $\mu = -J^{\bf A}H_{\bf A}$ is determined by $\Omega$ and the unknown functions $P(X,Y)$ and $Q(X,Y)$; however, one can obtain a simple determination of $\mu$ in terms of known quantities as follows. By Law \& Matsushita (2008) A1.6, the Ricci scalar curvature of the Walker geometry $(U,g,Z_{[\pi]})$ is $S = A_{UU} + B_{VV} + 2C_{UV} = \delta_{\bf A}\delta_{\bf B}{\cal W}^{\bf AB}$. The first term on the right-hand side of (3.43.7) expands as
$$\Omega^{-3}\delta^{\bf (A}[\Omega^2\delta^{\bf B)}\vartheta] = -2J^{\bf (A}\delta^{\bf B)}\vartheta + \Omega^{-1}\delta^{\bf A}\delta^{\bf B}\vartheta,\eqno(3.43.8)$$
whence $\delta_{\bf B}\left(\Omega^{-3}\delta^{\bf (A}[\Omega^2\delta^{\bf B)}\vartheta]\right) = J_{\bf B}\delta^{\bf A}\delta^{\bf B}\vartheta$. Hence, from (3.43.7) and using (A.10),
$$\eqalignno{S &= \delta_{\bf A}\delta_{\bf B}{\cal W}^{\bf AB} = \delta_{\bf A}\delta_{\bf B}\left({12\mu\Omega^{-3} + \hat S \over 12\tau^2}K^{\bf A}K^{\bf B}\right)\cr
& = {\mu \over \tau^2}K^{\bf A}K^{\bf B}\delta_{\bf A}\delta_{\bf B}\Omega^{-3}\cr
&= 6\mu\Omega^{-1}.&(3.43.9)\cr}$$
Substituting (3.43.9) into (3.43.7) yields (3.43.1).\bull
\vskip 24pt
The Einstein condition $\hat\Phi_{ABA'B'} = 0$ for AS$\alpha$-geometries is of special interest as then, by the GGST, the two conditions defining AS$\alpha$-geometries are equivalent to each other. Since $\hat\Phi_{ABA'B'} = 0$ entails, by the Bianchi identity, that $\hat S$ is constant, the Einstein condition for real AS$\alpha$-geometries is characterized by the conditions derived above together with $\hat\Phi_{AB1'1'} = 0$. In the hyperheavenly formalism, this final condition leads to a single PDE of the form $\delta_A\delta_B{\cal L} = 0$, where $\cal L$ is an expression in $\vartheta$, $\Omega$, and $\mu$ (equivalently $S$). Thus, $\cal L$ is affine in $U$ and $V$, with coefficients functions of $X$ and $Y$. This constraint is called the {\sl hyperheavenly equation} for $\vartheta$. We obtain this result as follows.

From (2.8), (3.15) and Law \& Matsushita (2008) (2.32),
$$\eqalignno{\hat\Phi_{AB1'1'} &= \hat\Phi_{ABA'B'}\hat\xi^{A'}\hat\xi^{B'}\cr
&= \Omega^{-3}\bigl[\Phi_{ABA'B'}\xi^{A'}\xi^{B'} + (\xi^{A'}\Upsilon_{AA'})(\xi^{B'}\Upsilon_{BB'}) - \xi^{B'}\xi^{A'}\nabla_{AA'}\Upsilon_{BB'}\bigr]\cr
&= \Omega^{-3}\bigl[A_{AB} + (\varsigma_A\omega)(\varsigma_B\omega) - \xi^{B'}\varsigma_A\nabla_{BB'}\omega\bigr]\cr
&= \Omega^{-3}\bigl[A_{AB} + 2\Omega^{-2}(\varsigma_A\Omega)(\varsigma_B\Omega) - \Omega^{-1}\xi^{B'}\varsigma_A\nabla_{BB'}\Omega\bigr].&(3.44)\cr}$$
In terms of the Walker spin frames, $\nabla_b\Omega = -(\delta_B\Omega)\xi_{B'} + (\varsigma_B\Omega)\pi_{B'}$, so the third term on the right-hand side of (3.44) may be expressed in terms of spin coefficients using Law (2008), (2.9), to obtain
$$\xi^{B'}\varsigma_A\nabla_{BB'}\Omega = \varsigma_A\varsigma_B\Omega + (\tilde\kappa'\alpha_A + \tilde\sigma'\beta_A)\delta_B\Omega + (\tilde\gamma\alpha_A -  \tilde\alpha\beta_A)\varsigma_B\Omega.\eqno(3.45)$$
We must now evaluate
$$\hat\Phi_{\bf AB1'1'} = \hat\Phi_{AB1'1'}\hat\epsilon_{\bf  A}{}^A\hat\epsilon_{\bf B}{}^B  = \Omega^{-1}\hat\Phi_{AB1'1'}\epsilon_{\bf A}{}^A\epsilon_{\bf B}{}^B,\eqno(3.46)$$
and do so by evaluating the terms in (3.44--45). From Law (2008), (5.11), one observes that for oriented Walker coordinates
$$A_{\bf AB} = \delta_{\bf(A}Z_{\bf B)},\qquad\hbox{where } Z_A := -(\tilde\sigma'\beta_A + \tilde\kappa'\alpha_A),\eqno(3.47)$$
is defined with respect to the corresponding (fixed) Walker spin frames (i.e., $Z_A$ is not well-defined under transformations between different sets of Walker spin frames). We now restrict to conformal oriented Walker coordinates $(U,V,X,Y)$ for $(M,h,[\pi^{A'}])$ so that (3.33-34), (A.10), and (3.40.1) all hold.

Using (3.40.1), one finds that, with 
$$\partial_{\bf B} := (-\partial_Y,\partial_X),\eqno(3.48)$$
$$\varsigma_A\varsigma_B\Omega = -{\Omega^2 \over 2}\bigl[-4\Omega^{-3}(\varsigma_A\Omega)(\varsigma_B\Omega) + \varsigma_A(\epsilon_B{}^{\bf F})\epsilon_{\bf FC}{\cal W}^{\bf CD}J_{\bf D} + \epsilon_B{}^{\bf F}\epsilon_{\bf FC}\epsilon_A{}^{\bf K}(\partial_{\bf K}{\cal W}^{\bf CD} + {1 \over 2}\epsilon_{\bf KP}{\cal W}^{\bf PQ}\delta_{\bf Q}{\cal W}^{\bf CD})J_{\bf D}\bigr].\eqno(3.49)$$
Writing $(\tilde\kappa'\alpha_A + \tilde\sigma'\beta_A)\delta_B\Omega = -Z_A\delta_B\Omega = \Omega^2Z_AJ_B$, then upon substituting this expression, (3.45), (3.47) and (3.49) into (3.44/46) and symmetrizing over {\bf A} and {\bf B}, one obtains
$$\eqalignno{\hat\Phi_{\bf AB1'1'} &= \Omega^{-4}\Bigl[\delta_{\bf (A}Z_{\bf B)} - \Omega Z_{\bf (A}J_{\bf B)} - {\Omega \over 2}\epsilon_{\bf C(B}\partial_{\bf A)}{\cal W}^{\bf CD}J_{\bf D}&(3.50)\cr
&\qquad + {\Omega \over 2}\bigl(\epsilon_{\bf (A}{}^A\epsilon_{\bf B)}{}^B\varsigma_A(\epsilon_B{}^{\bf F})\epsilon_{\bf FC}{\cal W}^{\bf CD}J_{\bf D} - {1 \over 2}\epsilon_{\bf C(B}\epsilon_{\bf A)P}{\cal W}^{\bf PQ}\delta_{\bf Q}{\cal W}^{\bf CD}J_{\bf D} + (\tilde\gamma\alpha_{\bf (A} - \tilde\alpha\beta_{\bf A})\epsilon_{\bf B)C}{\cal W}^{\bf CD}J_{\bf D}\bigr)\Bigr].\cr}$$
One computes:
$$-\epsilon_{\bf C(B}\partial_{\bf A)}{\cal W}^{\bf CD}J_{\bf D} -  J^{\bf C}\partial_{\bf C}{\cal W}_{\bf AB} = -\hat Z_{\bf (A}J_{\bf B)},\eqno(3.51)$$
where
$$\hat Z_{\bf B} := \pmatrix{C_Y-B_X\cr C_X-A_Y\cr}_{\bf B} = -\partial^{\bf C}{\cal W}_{\bf CB}.\eqno(3.52)$$
Writing $Z_{\bf A} =: (1/2)\hat Z_{\bf A} + Y_{\bf A}$, then
$$4Y_{\bf B} = \pmatrix{AB_U + C(B_V-C_U) - BC_V\cr -AC_U + C(A_U-C_V) + BA_V\cr}_{\bf B} = {\cal W}^{\bf CD}\delta_{\bf D}{\cal W}_{\bf BC}.\eqno(3.53)$$
Upon substituting (3.51), the first three terms inside the square brackets on the right-hand side of (3.50) may be rewritten in terms of $\hat Z_{\bf B}$ and $Y_{\bf B}$ as follows:
$$\Omega^{-1}\delta_{\bf (A}\Omega Y_{\bf B)} + {1 \over 2}\Omega^{-2}\delta_{\bf (A}\Omega^2\hat Z_{\bf B)} + {\Omega \over 2}J^{\bf C}\partial_{\bf C}{\cal W}_{\bf AB}.$$
Substituting for ${\cal W}_{\bf AB}$ using (3.42.1), but with the $L^{\bf B}$ of (3.43.6) replacing $F^{\bf B}$, and noting that the second summand of (3.42.1) is independent of $X$ and $Y$, one further obtains
$$\Omega^{-1}\delta_{\bf (A}\Omega Y_{\bf B)} + {\Omega^{-2} \over 2}[\delta_{\bf (A}\Omega^2\hat Z_{\bf B)} + J^{\bf C}\partial_{\bf C}\delta_{\bf (A}L_{\bf B)}].$$
Substituting now from (3.52) for $\hat Z_{\bf B}$, and noting that $\partial_{\bf B}$ and $\delta_{\bf B}$ are just partial derivatives and so commute, yields
$$\Omega^{-1}\delta_{\bf (A}\Omega Y_{\bf B)} - {\Omega^{-2} \over 2}\delta_{\bf (A}[\Omega^2\partial^{\bf C}{\cal W}_{\bf B)C} - J^{\bf C}\partial_{\bf |C|}L_{\bf B)}].\eqno(3.54)$$
Thus, (3.50) may be written
$$\eqalignno{\hat\Phi_{\bf AB1'1'} &= \Omega^{-4}\Bigl[- {\Omega^{-2} \over 2}\delta_{\bf (A}(\Omega^2\partial^{\bf C}{\cal W}_{\bf B)C} - J^{\bf C}\partial_{\bf |C|}L_{\bf B)}) + \Omega^{-1}\delta_{\bf (A}\Omega Y_{\bf B)}&(3.55)\cr
&\qquad + {\Omega \over 2}\bigl(\epsilon_{\bf (A}{}^A\epsilon_{\bf B)}{}^B\varsigma_A(\epsilon_B{}^{\bf F})\epsilon_{\bf FC}{\cal W}^{\bf CD}J_{\bf D} - {1 \over 2}\epsilon_{\bf C(B}\epsilon_{\bf A)P}{\cal W}^{\bf PQ}\delta_{\bf Q}{\cal W}^{\bf CD}J_{\bf D} + (\tilde\gamma\alpha_{\bf (A} - \tilde\alpha\beta_{\bf A})\epsilon_{\bf B)C}{\cal W}^{\bf CD}J_{\bf D}\bigr)\Bigr].\cr}$$
The first and third summands in the second line of (3.55) can be evaluated using Law (2008), (2.9) and (5.6). One finds
$$\epsilon_{\bf (A}{}^A\epsilon_{\bf B)}{}^B\varsigma_A(\epsilon_B{}^{\bf F})\epsilon_{\bf FC}{\cal W}^{\bf CD}J_{\bf D} + (\tilde\gamma\alpha_{\bf (A} - \tilde\alpha\beta_{\bf A})\epsilon_{\bf B)C}{\cal W}^{\bf CD}J_{\bf D} = -{1 \over 2}J_{\bf C}{\cal W}^{\bf CD}\delta_{\bf D}{\cal W}_{\bf AB}.\eqno(3.56)$$
For the remaining summand in the second line of (3.55), observe that
$$\eqalign{-{1 \over 2}\epsilon_{\bf CB}\epsilon_{\bf AP}{\cal W}^{\bf PQ}\delta_{\bf Q}{\cal W}^{\bf CD}J_{\bf D} &= {1 \over 2}\epsilon_{\bf AP}{\cal W}^{\bf PQ}\delta_{\bf Q}{\cal W}_{\bf BD}J^{\bf D}\cr
&= {\cal W}^{\bf PQ}\delta_{\bf Q}{\cal W}_{\bf B[A}J_{\bf P]}\cr
&= {1 \over 2}({\cal W}^{\bf PQ}\delta_{\bf Q}{\cal W}_{\bf BA}J_{\bf P} - {\cal W}^{\bf PQ}\delta_{\bf Q}{\cal W}_{\bf BP}J_{\bf A}).\cr}$$
Hence, upon substituting from (3.53), one finds
$$-{1 \over 2}\epsilon_{\bf C(B}\epsilon_{\bf A)P}{\cal W}^{\bf PQ}\delta_{\bf Q}{\cal W}^{\bf CD}J_{\bf D} = {1 \over 2}(J_{\bf P}{\cal W}^{\bf PQ}\delta_{\bf Q}{\cal W}_{\bf AB} - 4J_{\bf (A}Y_{\bf B)}).\eqno(3.57)$$
Hence, (3.56) cancels the first term on the right-hand side of (3.57), which entails that on the right-hand side of (3.55) the second line together with the last summand of the first line together take the form
$$\eqalignno{\Omega^{-1}\delta_{\bf (A}\Omega Y_{\bf B)} + {\Omega \over 2}[-2J_{\bf (A}Y_{\bf B)}] &= \delta_{\bf (A}Y_{\bf B)} - 2\Omega J_{\bf (A}Y_{\bf B)}\cr
&= \Omega^{-2}\delta_{\bf (A}\Omega^2Y_{\bf B)}\cr
&= {\Omega^{-2} \over 4}\delta_{\bf (A}\Omega^2{\cal W}^{\bf CD}\delta_{\bf |D|}{\cal W}_{\bf B)C},&(3.58)\cr}$$
upon substituting back with (3.53). Hence, substituting (3.58) into (3.55) yields
$$\hat\Phi_{\bf AB1'1'} = {\Omega^{-6} \over 4}\delta_{\bf (A}X_{\bf B)},\eqno(3.59)$$
where
$$X_{\bf B} := \Omega^2{\cal W}^{\bf CD}\delta_{\bf D}{\cal W}_{\bf BC} - 2(\Omega^2\partial^{\bf C}{\cal W}_{\bf BC} - J^{\bf C}\partial_{\bf C}L_{\bf B}).\eqno(3.60)$$
With slight differences, equations (3.59--60) are of the same form as Finley and Pleba\'nski (1976) (A4--5) when one formally identifies their $\phi$ with our $\Omega^{-1}$. We may therefore follow their argument from this point. First note that by (3.43.3) one can add a term $\lambda T_{\bf B}$, with $\lambda$ independent of $U$ and $V$, to $X_{\bf B}$ without altering (3.59), with the aim of choosing $\lambda$ so that $X_{\bf B}$ then equals $\delta_{\bf B}{\cal L}$ for some function $\cal L$.

Considering first the bracketed summands in (3.60), substituting (3.42.1) for ${\cal W}_{\bf BC}$, again with the $L_{\bf B}$ of (3.43.6) replacing $F_{\bf B}$, yields
$$\eqalign{-\Omega^2\partial^{\bf C}{\cal W}_{\bf BC} + J^{\bf C}\partial_{\bf C}L_{\bf B} &=  -\Omega^2\partial^{\bf C}\left(\Omega^{-3}\delta_{\bf (B}L_{\bf  C)}\right) -  J_{\bf C}\partial^{\bf  C}L_{\bf B}\cr
&= -{\Omega^{-1} \over  2}\left(\delta_{\bf B}\partial^{\bf C}L_{\bf C} +  \partial^{\bf C}\delta_{\bf C}L_{\bf B}\right) - J_{\bf C}\partial^{\bf C}L_{\bf B}\cr
&= -\delta_{\bf  B}\left(\Omega^{-1}\partial^{\bf C}L_{\bf C}\right) + {\Omega^{-1} \over 2}\left(\delta_{\bf B}\partial^{\bf C}L_{\bf C} - \delta_{\bf C}\partial^{\bf C}L_{\bf B}\right) + \left(J_{\bf B}\partial^{\bf C}L_{\bf C} - J_{\bf C}\partial^{\bf C}L_{\bf B}\right)\cr
&= -\delta_{\bf  B}\left(\Omega^{-1}\partial^{\bf C}L_{\bf C}\right) + {1 \over 2}\epsilon_{\bf BC}\left(\Omega^{-1}\delta_{\bf  D}\partial^{\bf C}L^{\bf D} + 2J_{\bf D}\partial^{\bf C}L^{\bf D}\right)\cr
&= -\delta_{\bf  B}\left(\Omega^{-1}\partial^{\bf C}L_{\bf C}\right) - {1 \over 2}\partial_{\bf B}(\Omega^{-1}\delta_{\bf D}L^{\bf D} + 2 J_{\bf D}L^{\bf D})\cr
&= -\delta_{\bf  B}\left(\Omega^{-1}\partial^{\bf C}L_{\bf C}\right) - {1 \over 2}\partial_{\bf B}\bigl(\Omega\delta_{\bf D}(\Omega^{-2}L^{\bf D})\bigr)\cr
&= -\delta_{\bf  B}\left(\Omega^{-1}\partial^{\bf C}L_{\bf C}\right) - {(\partial_{\bf B}\mu)T^{\bf D}K_{\bf D} \over \tau},\cr}$$
where the last line follows from (3.43.5). If one now chooses
$$\lambda := -{K^{\bf D}\partial_{\bf D}\mu \over 2\tau},\eqno(3.61)$$
and observes that
$$\eqalign{\delta_{\bf B}\left[{(K^{\bf D}T_{\bf D})(T^{\bf C}\partial_{\bf C}\mu) \over  2\tau}\right] &= {K^{\bf D}\epsilon_{\bf BD}(T^{\bf C}\partial_{\bf C}\mu) + (K^{\bf D}T_{\bf D})\epsilon_{\bf B}{}^{\bf C}\partial_{\bf C}\mu \over 2\tau}\qquad\hbox{by (3.43.3)}\cr
&=  {K^{\bf D}(T_{\bf D}\partial_{\bf B}\mu  - T_{\bf B}\partial_{\bf D}\mu)  + K^{\bf D}T_{\bf D}\partial_{\bf B}\mu \over 2\tau}\cr
&= {2K^{\bf D}T_{\bf D}\partial_{\bf B}\mu - T_{\bf B}K^{\bf D}\partial_{\bf D}\mu \over 2\tau},\cr}$$
then
$$-(\Omega^2\partial^{\bf C}{\cal W}_{\bf BC} - J^{\bf C}\partial_{\bf C}L_{\bf B}) + \lambda T_{\bf B} = \delta_{\bf B}\left[-\Omega^{-1}\partial^{\bf  C}L_{\bf C} + {K^{\bf D}T_{\bf D}T^{\bf C}\partial_{\bf C}\mu \over 2\tau}\right].\eqno(3.62)$$
To treat the first summand on the right-hand side of (3.60), observe that
$$\Omega^2{\cal W}_{\bf BC}{\cal W}^C{}_{\bf D} = \epsilon_{\bf BD}{\cal W},\qquad\hbox{where } {\cal W} := {1 \over 2}\Omega^2{\cal W}_{\bf BC}{\cal W}^{\bf BC}.\eqno(3.63)$$
Hence,
$$\delta_{\bf B}{\cal W} = -\delta^{\bf D}\epsilon_{\bf BD}{\cal W} = \delta_{\bf D}(\Omega^2{\cal W}_{\bf BC}{\cal W}^{\bf CD})$$
and
$$\Omega^2{\cal W}^{\bf CD}\delta_{\bf D}{\cal W}_{\bf BC} = \delta_{\bf B}{\cal W} - {\cal W}_{\bf BC}\delta_{\bf D}(\Omega^2{\cal W}^{\bf CD}).\eqno(3.64)$$
Using (3.43.7), one computes
$$\delta_{\bf D}(\Omega^2{\cal W}^{\bf CD}) = 2\Omega^3J^{\bf C}J_{\bf D}\delta^{\bf D}\vartheta  + {\mu K^{\bf C} \over \tau} - {\hat S\Omega^3K^{\bf C} \over 6\tau}.\eqno(3.65)$$
Expansion of the following first and third expressions confirms the identities:
$$\delta_{\bf (A}[\Omega^2\delta_{\bf B)}\vartheta]  = \Omega^2\delta_{\bf A}\delta_{\bf B}\vartheta  - 2\Omega^3J_{\bf (A}\delta_{\bf B)}\vartheta =  \Omega\delta_{\bf A}\delta_{\bf B}(\Omega\vartheta) - 2\Omega^4J_{\bf A}J_{\bf B}\vartheta.$$
Substituting into (3.43.7) yields
$${\cal W}_{\bf AB} = \Omega^{-2}\delta_{\bf A}\delta_{\bf B}(\Omega\vartheta) - 2\Omega J_{\bf A}J_{\bf B}\vartheta +{12\mu\Omega^{-3} + \hat S \over 12\tau^2}K_{\bf A}K_{\bf  B}.\eqno(3.66)$$
The second summand in (3.64) can therefore be evaluated by multiplying together (3.65--66). Upon doing so, there is a single term not involving $\mu$ or $\hat S$:
$$(2J^{\bf C}J_{\bf D}\delta^{\bf D}\vartheta)\bigl(\Omega\delta_{\bf B}\delta_{\bf C}(\Omega\vartheta)\bigr) = -\delta_{\bf B}\left[\bigl(J^{\bf C}\delta_{\bf C}(\Omega\vartheta)\bigr)^2\right].\eqno(3.67)$$
The terms involving $\mu$ are
$$\displaylines{{2\mu \over \tau}(J^{\bf D}\delta_{\bf D}\vartheta)K_{\bf B} - 2\mu\Omega\vartheta J_{\bf B} + {\mu\Omega^{-2} \over \tau}\delta_{\bf  B}\bigl(K^{\bf C}\delta_{\bf C}(\Omega\vartheta)\bigr)\hfill\cr
\hfill = {2\mu \over \tau}(J^{\bf D}\delta_{\bf D}\vartheta )K_{\bf B} - {\mu \over \tau}(K^{\bf C}\delta_{\bf C}\vartheta)J_{\bf B} - \mu\delta_{\bf B}\vartheta +  {\mu\Omega^{-1} \over \tau}\delta_{\bf B}(K^{\bf C}\delta_{\bf C}\vartheta),\cr}$$
upon expanding out the last summand of the left-hand side and rearranging terms. The right-hand side can be written
$$-{2\mu\over \tau}(K^{\bf D}J_{\bf B} - J^{\bf D}K_{\bf B})\delta_{\bf D}\vartheta + {\mu \over \tau}J_{\bf B}K^{\bf D}\delta_{\bf D}\vartheta - \mu\delta_{\bf B}\vartheta + {\mu\Omega^{-1} \over \tau}\delta_{\bf B}(K^{\bf D}\delta_{\bf D}\vartheta),$$
which, upon utilizing (A10), becomes
$$-3\mu\delta_{\bf B}\vartheta + {\mu \over \tau}J_{\bf B}K^{\bf D}\delta_{\bf D}\vartheta + {\mu\Omega^{-1} \over \tau}\delta_{\bf B}(K^{\bf C}\delta_{\bf C}\vartheta) = {\mu \over \tau}\delta_{\bf B}\bigl(\Omega^{-4}K^{\bf D}\delta_{\bf D}(\Omega^3\vartheta)\bigr),\eqno(3.68)$$
as is easily checked by expanding out the right-hand side of (3.68).

Finally, the terms involving $\hat S$ are:
$${\hat S \over 6}\left[{\Omega^3 \over \tau}(J^{\bf D}\delta_{\bf D}\vartheta)K_{\bf B} - {\Omega \over \tau}\delta_{\bf B}\bigl(K^{\bf C}\delta_{\bf C}(\Omega\vartheta)\bigr) + 2\Omega^4\vartheta J_{\bf B}\right]$$
which, expanding out the middle term, yields
$${\hat S \over 6}\left[{\Omega^3 \over \tau}(J^{\bf D}K_{\bf B} + J_{\bf B}K^{\bf D})\delta_{\bf D}\vartheta - {\Omega^2 \over \tau}\delta_{\bf B}(K^{\bf D}\delta_{\bf D}\vartheta) + \Omega^3\delta_{\bf B}\vartheta\right],$$
which in turn can be written
$${\hat S \over 6}\left[{\Omega^3 \over \tau}(J^{\bf D}K_{\bf B} - K^{\bf D}J_{\bf B})\delta_{\bf D}\vartheta +{2\Omega^3 \over \tau}J_{\bf B}K^{\bf D}\delta_{\bf D}\vartheta - {\Omega^2 \over \tau}\delta_{\bf B}(K^{\bf D}\delta_{\bf D}\vartheta) + \Omega^3\delta_{\bf B}\vartheta\right].$$
Upon using (A10), one finds the last expression simplifies to
$$-{\hat S \over 6\tau}\delta_{\bf B}(\Omega^2K^{\bf D}\delta_{\bf D}\vartheta).\eqno(3.69)$$

Thus, (3.59) remains valid with (3.60) replaced by
$$X_{\bf B} := \Omega^2{\cal W}^{\bf CD}\delta_{\bf D}{\cal W}_{\bf BC} + 2(-\Omega^2\partial^{\bf C}{\cal W}_{\bf BC} + J^{\bf C}\partial_{\bf C}L_{\bf B} - {K^{\bf D}\partial_{\bf D}\mu \over 2\tau}T_{\bf B}) = \delta_{\bf B}{\cal L}\eqno(3.70)$$
where, by (3.62), (3.64) and (3.67--69),
$${\cal L} = -2\Omega^{-1}\partial^{\bf D}L_{\bf D} + {(K^{\bf D}T_{\bf D})(T^{\bf C}\partial_{\bf C}\mu) \over \tau} + {\cal W} -\left[-\bigl(J^{\bf C}\delta_{\bf C}(\Omega\vartheta)\bigr)^2 + {\mu \over \tau}\Omega^{-4}K^{\bf D}\delta_{\bf D}(\Omega^3\vartheta) - {\hat S \over 6\tau}\Omega^2K^{\bf D}\delta_{\bf D}\vartheta\right].$$
Defining 
$$w := \Omega^{-1} = MU + NV,\qquad\hbox{ then }\qquad K^{\bf D}\delta_{\bf D} = -{d \over dw},\eqno(3.71)$$
and substituting from (3.43.6) for $L_{\bf C}$ in favour of $\vartheta$, one can write $\cal L$ as
$${\cal L} = {\cal W} + \bigl(J^{\bf D}\delta_{\bf D}(\Omega\vartheta)\bigr)^2 + 2\Omega\partial_{\bf D}\delta^{\bf D}\vartheta - {\mu\Omega^{-4} \over \tau}{d(\Omega^3\vartheta) \over dw} + {T^{\bf D}K_{\bf D} \over \tau^2}\left[2\Omega^{-1}K^{\bf C} - \tau T^{\bf C}\right]\partial_{\bf C}\mu  + {\hat S\Omega^2 \over 6\tau}{d\vartheta \over dw}.\eqno(3.72)$$
Hence, in terms of the conformal oriented Walker coordinates $(U,V,X,Y)$ for which $\Omega^{-1} = MU + NV =: w$, ($M$ and $N$ constant), from (3.59) and (3.70)
$$\hat\Phi_{\bf AB1'1'} = {\Omega^{-6} \over 4}\delta_{\bf A}\delta_{\bf B}{\cal L} = 0 \iff \delta_{\bf A}\delta_{\bf B}{\cal L} = 0,\qquad\hbox{i.e., $\cal L$ is affine in $U$ and $V$:}\qquad {\cal L} = T^{\bf D}\eta_{\bf D} + \Xi,\eqno(3.73)$$
for constants $\eta_{\bf C}$ and $\Xi$, which is our version of the {\sl hyperheavenly equation}, cf. Finley \& Pleba\'nski (1976), (3.14).
 
Consequences of various conditions imposed on the Ricci curvature of a real AS$\alpha$-geometry $(M,h,[\pi^{A'}])$ derived in this section may be summarized as follows:\hfil\break
1) $(M,h,[\pi^{A'}])$ may be characterized locally by four arbitrary functions $a$, $b$, $c$, and $\Omega$ of local coordinates $(u,v,x,y)$ with the metric given with respect to these coordinates by (3.3), with $W$ as in (3.2);\hfil\break
2) $(M,h,[\pi^{A'}])$, with $[\pi^{A'}]$ a RPS, may be characterized locally as in (1) but with $\Omega$ an affine function of $u$ and $v$; one may specialize this local characterization to local coordinates $(U,V,X,Y)$ for which $\Omega = MU+NV$, with $M$ and $N$ constant;\hfil\break
3) $(M,h,[\pi^{A'}])$, with constant scalar curvature and with $[\pi^{A'}]$ a multiple RPS, may be characterized locally by three functions $\vartheta$, $\Omega$, and $S$ of local coordinates $(U,V,X,Y)$ such that $\vartheta$ is arbitrary, $\Omega =MU+NV$, $M$ and $N$ constant, $S\Omega$ is a function of $(X,Y)$ only, and the metric with respect to $(U,V,X,Y)$ is given by (3.3) with $W$ given by (3.43.1); the condition on $S\Omega$ ensures that the (locally defined) Walker metric $g := \Omega^{-2}h$ has scalar curvature $S$;\hfil\break
4) $(M,h,[\pi^{A'}])$ Einstein may be locally characterized as in (3), with $\vartheta$ subject to a single constraint equation, the hyperheavenly equation (3.72--73).

One can obviously adapt applications of the hyperheavenly equation in complex general relativity to the study of four-dimensional neutral geometry. In particular, various solutions of the hyperheavenly equation, e.g., Pleba\'nski \& Torres del Castillo (1982), will be of interest. More recently, Chudecki \& Przanowski (2008a) obtained an explicit neutral metric for a real AS$\alpha$-geometry as a solution of the hyperheavenly equation. 
\vskip 24pt
\noindent {\section 4. Null Geometry}
\vskip 12pt
$(M,h,[\pi^{A'}])$ again denotes a real AS$\alpha$-geometry and we shall continue to denote its curvature quantities etc., by hatted symbols. When convenient, one may suppose that any point $p \in M$ has a neighbourhood $U$ such that $h = \Omega^2g$ on $U$, for some positive function $\Omega$, with $(U,g,Z_{[\pi]})$ a Walker geometry. 

Suppose $(M,h,[\pi^{A'}])$ is  {\sl not} Walker. As already noted, (1.3) entails that $S^b = \pi_{A'}\hat\nabla^b\pi^{A'} = \omega^B\pi^{B'} \not= 0$, for some nonzero spinor $\omega^B$, where $\pi^{A'}$ is an LSR of $[\pi^{A'}]$. While the geometry only determines $S^b$ up to scale, the null distribution ${\cal D} := \langle S^b \rangle_{\bf R}$ is well defined. Indeed, the geometry defines the nested null distributions ${\cal D} \leq Z_{[\pi]} \leq {\cal H}$, where ${\cal H} = {\cal D}^\perp$, (these distributions are null of types I, II, and III, respectively, in the sense of Law 2008). The geometry of these nested distributions was considered in Law (2008) for $\alpha$-geometries. By Law (2008) (6.2.37), the condition for $\cal D$ to be auto-parallel, i.e., for any local section $q^a$ of $\cal D$, $q^b\hat\nabla_bq^a \propto q^a$, is $\hat\Phi_{ABA'B'}\omega^A\omega^B\pi^{A'}\pi^{B'} = 0$. If one supposes that $[\pi^{A'}]$ is a RPS, then, by Law (2008) (6.2.18), an $\alpha$-geometry is in fact an AS$\alpha$-geometry. Furthermore, by the proof of (6.3.14) in Law (2008), in an AS$\alpha$-geometry, $\hat\Phi_{ABA'B'}\omega^B\pi^{A'}\pi^{B'} = 0$ is the condition for $\cal H$ to be integrable. 

Hence, an (AS)$\alpha$-geometry $(M,h,[\pi^{A'}])$ for which $[\pi^{A'}]$ is a RPS, which (consistent with Law 2008, (6.3.16)) we call a {\sl Ricci-aligned} (AS)$\alpha$-geometry, has nested integrable null distributions ${\cal D} \leq Z_{[\pi]} \leq {\cal H}$, with $\cal D$ auto-parallel. The integral curves of $\cal D$, suitably parametrized, are the null geodesic generators of the null hypersurfaces which are the integral hypersurfaces of $\cal H$. These null hypersurfaces are also foliated by the $\alpha$-surfaces of $Z_{[\pi]}$. Frobenius coordinates for these nested distributions were described in Law (2008) (6.3.16--18) and provide an alternative to conformal Walker coordinates. 

Note that as distributions, each of $\cal D$, $Z_{[\pi]}$, and $\cal H$ is of course defined and integrable in $(U,g)$. Moreover each retains its null character as these are conformally invariant. As $\cal D$ is null, the condition of being auto-parallel is conformally invariant too, so $\cal D$ is auto-parallel with respect to $g$.

By Law (2008) (6.2.29), for any spin frames $\{o^A,\iota^A\}$ and $\{o^{A'},\iota^{A'}\}$ for which $o^{A'}$ is an LSR of $[\pi^{A'}]$, $\omega_A = \tilde\tau o_A - \tilde\rho\iota_A$. Hence, with respect to the spin frames $\{\hat\alpha^A,\hat\beta^A\}$ and $\{\hat\xi^{A'},\hat\pi^{A'}\}$ of (3.15) for $(M,h,[\pi^{A'}])$ associated to conformal oriented Walker coordinates $(u,v,x,y)$, by (3.21--22) 
$$\eqalignno{\omega_A &= \hat{\tilde\tau}\hat\alpha_A - \hat{\tilde\rho}\hat\beta_A\cr
&= \Omega^{1/2}(\triangle\omega)\alpha_A - \Omega^{1/2}(D\omega)\beta_A\cr
&= \Omega^{-1/2}\left[(\triangle\Omega)\alpha_A - (D\Omega)\beta_A\right]\cr
&= \Omega^{-1/2}\delta_A\Omega.&(4.1)\cr}$$
The local conformal geometry also defines $\hat\nabla_a\Omega = \nabla_a\Omega$, i.e., the distribution $\langle \Upsilon^a \rangle_{\bf R}$. In terms of the conformally associated Walker spin frames, $\nabla_a\Omega = -(\delta_A\Omega)\xi_{A'} + (\varsigma_A\Omega)\pi_{A'}$. Note that $(M,h,[\pi^{A'}])$ is Walker iff $\omega_A$ vanishes, i.e., iff $\delta_A\Omega$ vanishes, i.e., in accordance with (3.1), iff $\pi^{A'}\nabla_{AA'}\Omega$ vanishes, in which case $\nabla_a\Omega = (\varsigma_A\Omega)\pi_{A'}$ and is null. More generally, for any Ricci-aligned $\alpha$-geometry,
$$\eqalign{(\nabla_a\Omega)(\nabla^a\Omega) &= -2(\delta^A\Omega)(\varsigma_A\Omega)\cr
&= -{\cal W}^{\bf CD}(\delta_{\bf C}\Omega)(\delta_{\bf D}\Omega),\qquad\hbox{by (3.40.1).}\cr}$$
One does not expect, generically, $\nabla_a\Omega$ to be null. But, for any AS$\alpha$-geometry
$$(\nabla_a\Omega)S^a \propto \bigl(-(\delta_A\Omega)\xi_{A'} + (\varsigma_A\Omega)\pi_{A'}\bigr)\delta^A\Omega\pi^{A'} = 0.$$
Thus, if $(M,h,[\pi^{A'}])$ is not Walker, $S^a$ does not vanish, $\cal H$ is well defined, and $\nabla^a\Omega \in {\cal H}$, which thereby relates this ingredient of the local conformal geometry to the geometry of the nested null distributions ${\cal D} \leq Z_{[\pi]} \leq {\cal H}$. By (4.1), the null distribution $\cal D$ is aligned with the local conformal Walker geometry, i.e., $\omega_A$ is proportional to $\alpha_A$ or $\beta_A$, iff $D\Omega =0$ or $\triangle\Omega = 0$, respectively, i.e., $\Omega$ is independent of $U$ or $V$, respectively.

If $U$ is a domain of $M$ on which there is an LSR $\pi^{A'}$ and $\Omega$ and $\chi$ are local conformal factors whose domains intersect in $U$ then, by (4.1), $\delta_A(\Omega^{1/2}) = \delta_A(\chi^{1/2})$, i.e., $\Omega^{1/2}$ and $\chi^{1/2}$ differ only by a constant on a given $\alpha$-surface within $U$. 
\vskip 24pt
\noindent {\section  Appendix}
\vskip 24pt
Let $(u,v,x,y)$ and $(p,q,r,s)$ be two sets of overlapping oriented Walker coordinates for a Walker geometry $(M,g)$. Suppose the metric takes components with respect to $(u,v,x,y)$ as in (3.2), and with respect to $(p,q,r,s)$ of the same form but with $\check W$ in place of $W$ and $\check a$, $\check b$, and $\check c$ in place of $a$, $b$, and $c$ respectively. For each set of coordinates one can construct the corresponding Walker null tetrads and spin frames as in (3.6). The  notation employed in (3.6) will denote the Walker null tetrad and spin frames for $(u,v,x,y)$. The Walker null tetrad and spin frames for $(p,q,r,s)$ will be distinguished by the use of the `check' mark over the relevant symbol. For notational convenience, the pair of Walker spin frames will also be denoted here by $\epsilon_{\bf A}{}^A$ and $\epsilon_{\bf A'}{}^{A'}$. Let $\{\partial_1,\ldots,\partial_4\}$ denote the coordinate basis for $(u,v,x,y)$ and $\{\flat_1,\ldots,\flat_4\}$ that for $(p,q,r,s)$. Then
$$\flat_j = \sum_{i=1}^4\,\partial_iJ^i{}_j,\qquad\hbox{where}\qquad (J^i{}_j) := {\partial(u,v,x,y) \over \partial(p,q,r,s)} = \pmatrix{D&E\cr {\bf 0}_2&{^\tau\! D}^{-1}\cr},\eqno({\rm A}1)$$
with $E$, $D \in {\bf R}(2)$ and $\det(D) > 0$, see Law \& Matsushita (2008) (1.2). Our aim in this Appendix is to record the relationships between the Walker spin frames for the two sets of oriented Walker coordinates and some simple observations. From (A1) and (3.6), one can express the two Walker null tetrads in terms of each other and then deduce the relationship between the Walker spin frames, obtaining:
$$\vcenter{\openup1\jot \halign{$\hfil#$&&${}#\hfil$&\qquad$\hfil#$\cr
\check\epsilon_{\bf B}{}^A &= \epsilon_{\bf A}{}^A\Lambda^{\bf A}{}_{\bf B} & \Lambda &:= \left(\Lambda^{\bf A}{}_{\bf B}\right) = \chi^{-1}D \in {\bf SL(2;R)}&\chi := \pm\sqrt{\det(D)}\cr
\check\epsilon_{\bf B'}{}^{A'} &= \epsilon_{\bf A'}{}^{A'}\tilde\Lambda^{\bf A'}{}_{\bf B'} & \tilde\Lambda &:= \left(\tilde\Lambda^{\bf A'}{}_{\bf B'}\right) = \pmatrix{\chi&\chi^{-1}\mu\cr 0&\chi^{-1}\cr} \in {\bf SL(2;R)}\cr}}\eqno({\rm A}2)$$
where
$$\eqalignno{\mu &:= {(D^1{}_1D^2{}_2 + D^2{}_1D^1{}_2)c - D^1{}_1D^1{}_2b - D^2{}_2D^2{}_1a \over 2\chi^2} - {\chi^2\check c \over 2} + D^1{}_1E^2{}_1 - D^2{}_1E^1{}_1\cr
&= {\chi^2\check c \over 2} - {(D^1{}_1D^2{}_2 + D^2{}_1D^1{}_2)c - D^1{}_1D^1{}_2b - D^2{}_2D^2{}_1a \over 2\chi^2} +D^1{}_2E^2{}_2 - D^2{}_2E^1{}_2.&({\rm A}3)\cr}$$
The ambiguity in sign for $\chi$ corresponds to the ambiguity in sign of the LSR of $[\pi^{A'}]$ satisfying (3.5) (Law \& Matsushita 2008, (2.8)) and the ambiguity in overall sign for the Walker spin frames. Note that, by assumption,
$$\pmatrix{{\bf 0}_2&{\bf 1}_2\cr {\bf 1}_2&\check W\cr} = {^\tau\! J}.\pmatrix{{\bf 0}_2&{\bf 1}_2\cr {\bf 1}_2&W\cr}.J.\eqno({\rm A}4)$$
The nontrivial condition in (A4) is
$$\check W = {^\tau\! E}.{^\tau\! D}^{-1} + D^{-1}.E + D^{-1}.W.{^\tau\! D}^{-1}.\eqno({\rm A}5)$$
The equality in (A3) is a consequence of the equation in the off-diagonal terms of (A5).

The quantity $\delta_A := \pi^{A'}\nabla_{AA'}$ where $\pi^{A'}$ is any LSR for $[\pi^{A'}]$ is clearly determined, up to scale, by the (Walker) geometry. Fixing the LSR $\pi^{A'}$ to be that in (3.5), i.e., to be the element $\epsilon_{0'}{}^{A'}$ of the Walker spin frames, fixes $\delta_A$ up to sign. Moreover, from (A1--2), one confirms that $\check\delta^A = \pm \chi\delta_A$, as expected. On the other hand, defining, with respect to a given pair of Walker spin frames,
$$\varsigma_A := \xi^{A'}\nabla_{AA'} = \alpha_AD' - \beta_A\delta,\qquad\hbox{then}\qquad \check\varsigma_A = \chi^{-1}(\varsigma_A + \mu\delta_A).\eqno({\rm A}6)$$
Using Law (2008), (5.8) \& (5.10), one can show that, acting on functions,
$$\delta_{(A}\varsigma_{B)} = \varsigma_{(A}\delta_{B)} - [\varsigma_{(A}\pi^{B'}]\nabla_{B)B'} .\eqno({\rm A}7)$$
The hyperheavenly formalism, e.g., of Finley and Pleba\'nski (1976) and Boyer et al. (1980), might tempt one to define, for given Walker spin frames, $W^{AB} := W^{\bf AB}\epsilon_{\bf A}{}^A\epsilon_{\bf B}{}^B$. From (A2--3) and (A5), one computes
$$\check W^{AB} = \chi^{-2}\bigl((D.{^\tau\! E} + E.{^\tau\! D})^{\bf CD}\epsilon_{\bf C}{}^A\epsilon_{\bf D}{}^B + W^{AB}\bigr).\eqno({\rm A}8)$$
Hence, for $W^{AB}$ to be a meaningful spinorial object, at least with respect to Walker spin frames, one would require
$$D \in {\bf SL(2;R)}\qquad\hbox{and}\qquad D.{^\tau\! E}\hbox{  skew}.\eqno({\rm A}9)$$

When utilizing oriented Walker coordinates $(U,V,X,Y)$ as in (3.33--34), and their associated Walker spin frames, it will be convenient to employ the following notation. Noting (3.34--35), define 
$$\displaylines{J_A := \delta_A\Omega^{-1} = -\Omega^{-2}\delta_A\Omega; \quad K^A := -(N\alpha^A + M\beta^A); \qquad \tau := K^AJ_A = -2MN\cr
\hfill\hbox{whence}\qquad 2K^{[A}J^{B]} = -\tau\epsilon^{AB}.\hfill\llap({\rm A}10)\cr}$$
If $(P,Q,R,S)$ are another set of oriented Walker coordinates satisfying (3.33--34), then
$$\Omega = MU+NV = FP+GQ,\eqno({\rm A}11)$$
say. Consequently, $MdU + NdV = FdP + GdQ$, which is equivalent to, with $J := J_A\epsilon_{\bf A}{}^A$,
$${^\tau\! D}.J = \check J \hskip 1.25in {^\tau\! E}.J = {\bf 0},\eqno({\rm A}12)$$
the first equation of which is consistent with $\check\delta_A\Omega = \chi\delta_A\Omega$. $K^A$, however, does not have form-preserving transformation and is best thought of as defined with respect to a fixed set of such oriented Walker coordinates $(U,V,X,Y)$, and the associated Walker spin frames.

\vskip 24pt
\noindent {\section References}
\frenchspacing
\vskip 1pt
\hangindent=20pt \hangafter=1
\noindent Boyer, C. P., Finley III, J. D. \& Pleba\'nski, J. F. 1980 Complex General Relativity, $\cal H$ and ${\cal HH}$ Spaces-A Survey of One Approach, in {\sl General Relativity and Gravitation: One Hundred Years After the Birth of Albert Einstein}, Vol. 2, A. Held (ed.), Plenum Press, New York, NY, 241--281.
\vskip 1pt
\hangindent=20pt \hangafter=1
\noindent Chudecki, A. \& Przanowski, M. 2008a A simple example of type-$[{\rm N}] \otimes [{\rm N}]{\cal HH}$-spaces admitting twisting geodesic congruence. {\sl Classical Quantum Gravity} {\bf 25}, 055010 (13pp).
\vskip 1pt
\hangindent=20pt \hangafter=1
\noindent Chudecki, A. and Przanowski, M. 2008b From hyperheavenly spaces to Walker spaces and Osserman spaces. I. {\sl Classical Quantum Gravity}, {\bf 25} (145010) (18 pp).
\vskip 1pt
\hangindent=20pt \hangafter=1
\noindent Finley III, J. D. \& Pleba\'nski, J. F. 1976 The intrinsic spinorial structure of hyperheavens. {\sl Journal of Mathematical Physics} {\bf 17}, 2207--2214.
\vskip 1pt
\hangindent=20pt \hangafter=1
\noindent Law, P. R. 2006 Classification of the Weyl curvature spinors of neutral metrics in four dimensions. {\sl J. Geo. Phys.} {\bf 56}, 2093--2108.
\vskip 1pt
\hangindent=20pt \hangafter=1
\noindent Law, P. R. 2008 Spin Coefficients for Four-Dimensional Neutral Metrics, and Null Geometry. {\sl J. Geo. Phys.}, {\bf 59(8)}, 1087--1126. arXiv:0802.1761v2 [math.DG] (26 Aug 2009).
\vskip 1pt
\hangindent=20pt \hangafter=1
\noindent Law, P. R. \& Matsushita, Y. 2008 A Spinor Approach to Walker Geometry. {\sl Communications in Mathematical Physics} {\bf 282}, 577--623. arxiv:math/0612804v4 [math.DG] (7 Apr 2009).
\vskip 1pt
\hangindent=20pt \hangafter=1
\noindent Penrose, R. \& Rindler, W. 1984 {\sl Spinors and Space-Time, Vol. 1: Two-spinor calculus and relativistic fields\/}, Cambridge University Press, Cambridge.
\vskip 1pt
\hangindent=20pt \hangafter=1
\noindent Penrose, R. \& Rindler, W. 1986 {\sl Spinors and Space-Time, Vol. 2: Spinor and twistor methods in space-time geometry\/}, Cambridge University Press, Cambridge.
\vskip 1pt
\hangindent=20pt \hangafter=1
\noindent Pleba\'nski, J. F. \& Robinson, I. 1976 Left-Degenerate Vacuum Metrics. {\sl Physical Review Letters} {\bf 37}, 493--495.
\vskip 1pt
\hangindent=20pt \hangafter=1
\noindent Pleba\'nski, J. F. \& Robinson, I. 1977 The Complex Vacuum Metric with Minimally Degenerated Conformal Curvature; in {\sl Asymptotic Structure of Space-Time}, F. P. Esposito \& L. Witten (eds), Plenum Press, New York \& London, 361--406.
\vskip 1pt
\hangindent=20pt \hangafter=1
\noindent Pleba\'nski, J. F. \& R\'ozga, K. 1984 The optics of null strings. {\sl Journal of Mathematical Physics} {\bf 25}, 1930--1940.
\vskip 1pt
\hangindent=20pt \hangafter=1
\noindent Pleba\'nski, J. F. \& Torres del Castillo, G. F. 1982 ${\cal HH}$ spaces with an algebraically degenerate right side. {\sl Journal of Mathematical Physics} {\bf 23}, 1349--1352.
\vskip 1pt
\hangindent=20pt \hangafter=1
\noindent Walker, A. G. 1950 Canonical form for a Riemannian space with a parallel field of null planes. {\sl Quart. J. Math. Oxford(2)\/} {\bf 1}, 69--79.
\vskip 1pt

\bye